\documentclass[10pt,a4paper]{article}
\usepackage[]{graphicx,float,latexsym,times}
\usepackage{amsfonts,amstext,amsmath,amssymb,amsthm}

\RequirePackage[OT1]{fontenc}
\long\def\symbolfootnote[#1]#2{\begingroup%
	\def\thefootnote{\fnsymbol{footnote}}\footnote[#1]{#2}\endgroup}
\usepackage{bm}
 
\usepackage{xspace}
\usepackage{geometry}
\usepackage{array}
\usepackage{subfigure}
\usepackage{color}
\usepackage{epsfig}
\usepackage{hhline}
\usepackage{natbib}
\usepackage{bbm}
\usepackage[T1]{fontenc}
\usepackage[]{graphicx,float,latexsym,times}
\usepackage{amsfonts,amstext,amsmath,amssymb}
\usepackage[utf8]{inputenc}    % pour les accents 
\DeclareMathAccent{\widehat}{\mathord}{largesymbols}{"62}
\DeclareMathAccent{\widetilde}{\mathord}{largesymbols}{"65}
\def\pth#1{\left(#1\right)}

\newcolumntype{B}{>{\usefont{T1}{ptm}{m}{n}\fontsize{.3cm}{.3cm}\selectfont }c} % font Times-Roman
\def\ebo{\textrm{\mathversion{bold}$\mathbf{\beta}^0$\mathversion{normal}}}

\def\oo{\textrm{\mathversion{bold}$\mathbf{0}$\mathversion{normal}}}
\def\eb{\textrm{\mathversion{bold}$\mathbf{\beta}$\mathversion{normal}}}

\def\eo{\textrm{\mathversion{bold}$\mathbf{\omega}$\mathversion{normal}}} 
\def\ek{\textrm{\mathversion{bold}$\mathbf{\kappa}$\mathversion{normal}}} 
\def\evs{\textrm{\mathversion{bold}$\mathbf{\varsigma}$\mathversion{normal}}} 
\def\eS{\textrm{\mathversion{bold}$\mathbf{\varUpsilon}$\mathversion{normal}}}

\def\evP{\textrm{\mathversion{bold}$\mathbf{\varPhi}$\mathversion{normal}}}

\def\eE{\mathbb{E}}

\def\e1{1\!\!1}

\newtheorem{theorem}{Theorem} 
\def\ee1{\textrm{\mathversion{bold}$\mathbf{\varepsilon}$\mathversion{normal}}}

\def\oo{\textrm{\mathversion{bold}$\mathbf{0}$\mathversion{normal}}}

\newcommand{\eV}{{\bf V}}
\newcommand{\eD}{{\bf D}}

\newcommand{\bfQ}{{\bf Q}}
\newcommand{\ev}{{\bf v}}
\newcommand{\eW}{{\bf W}}
\newcommand{\eI}{{\bf I}}
\newcommand{\eh}{{\bf h}}

\def\eu{\mathbf{{u}}}
\def\es{\mathbf{{s}}}

\newcommand{\N}{\mathbb{N}}
\newcommand{\R}{\mathbb{R}}
\newcommand{\PP}{\mathbb{P}}
\newcommand{\eSS}{\mathbb{S}}

\def\eX{\mathbf{X}}
\def\ex{\mathbf{x}}

\newcommand{\Var}{\mathbb{V}\mbox{ar}\,}

\def\argmin{\mathop{\mathrm{arg\,min}}} 

\begin{document}
%--------------------------------------------------
%--------------------------------------------------
 
\title {{Right-censored models on massive data}}
 
%--------------------------------------------------
    \date{}
	\author{Gabriela CIUPERCA}
\setcounter{footnote}{1}
\footnotetext{\noindent \small{\indent Universite Claude Bernard Lyon 1, ICJ UMR 5208, 
		Bat.  Braconnier, 43, blvd du 11 novembre 1918, F - 69622 Villeurbanne Cedex, France, 
		\textit{Email address}: Gabriela.Ciuperca@univ-lyon1.fr}}
% \footnote{\noindent \small{CONTACT: Gabriela CIUPERCA,\\  Universit\'e de Lyon,Universit\'e Claude Bernard Lyon 1, CNRS, UMR 5208, Institut Camille Jordan, Bat.  Braconnier, M43, blvd du 11 novembre 1918, F - 69622 Villeurbanne Cedex, France\\
		%	\textit{Email address}: Gabriela.Ciuperca@univ-lyon1.fr }}
\maketitle
%%%%%%%%%%%%%%%
\begin{abstract}
 This article considers the automatic selection  problem of the relevant explanatory variables in a right-censored model on a massive database. We propose and study four aggregated censored adaptive LASSO estimators constructed by dividing the observations in such a way as to keep the consistency of the estimator of the survival curve. We show that these estimators have the same theoretical oracle properties as the one built on the full database. Moreover, by Monte Carlo simulations we obtain that their calculation time is smaller than that of the full database. The simulations confirm  also the theoretical properties. For optimal tuning parameter selection, we propose a BIC-type criterion.
 \end{abstract}
\noindent \textbf{Keywords}:   massive data, right-censoring, automatic selection, adaptive LASSO, BIC criterion.

%%%%%%%%%%%%%%%

%%%%%%%%%%%%%%%% 
\section{Introduction}
Statistical inference for a model may prove impossible in practice if the database is massive. This is all the more true and important for a censoring model, where the survival curve estimation must be carried out at the same time as automatic selection of the relevant variables.  Therefore, heavy numerical calculations must be undertaken which will be impossible to carry out due to the computer's limited memory. To the author's knowledge, very few works deal with censored models  when the sample size is extraordinarily large and when the number of explanatory variables is large, but smaller than the sample size. The very recent paper of \cite{Su.Yin.Zhang.Zhao.23} builds a penalized loss function by approximating the weighted least squares loss function by combining estimation results without penalization from all sets. Then the estimator calculated on $K$ groups of observations considered by \cite{Su.Yin.Zhang.Zhao.23} is totally different from the aggregated estimators proposed in the present paper. Moreover, their model error was considered to be of zero expectation. If the number of explanatory variables is large, then automatic selection of significant variables is required. For classical linear regression models, by penalising a loss function with an adaptive LASSO penalty, as initially proposed by \cite{Zou.06}, the obtained estimators satisfy the oracle properties, i.e. the sparsity (that is the non-zero coefficients are shrinked directly as 0 with a probability converging to 1) and asymptotic normality of the estimators of non-zero coefficients. These adaptive LASSO penalties have been considered for automatic selection in censored models for several types of loss functions.\\
\cite{Wang.Zhou.Li.13} developed a variable selection method for censored quantile regression based on the adaptive LASSO penalisation for  a quantile process weighted by weights to redistribute mass. A BIC-type criterion is also proposed for selecting the penalty tuning parameter. An extended BIC criterion was proposed by \cite{Lee.Park.Lee.23} for a quantile forward regression  for ultrahigh-dimensional survival model, while in \cite{Wang.Song.2011} a BIC criterion is proposed for the AFT model estimated by the adaptive LASSO method for the least squares (LS) loss function. A high-dimensional censored quantile regression with a divergent number of explanatory variables and where the number of non-zero coefficients may depend on the number of observations  is also considered by \cite{Fei.Zheng.Hong.Li.23}. For censored models, if the quantile index of errors is unknown, the adapted LASSO composite quantile method proposed and studied by \cite{Tang-Zhou-Wu.12} can be used. For LS loss function with adaptive LASSO penalty, \cite{Wang.Song.2011} proves the oracle properties of the corresponding estimator and proposes an algorithm for its calculation. Afterwards, \cite{Ciuperca.24} generalizes the estimation by censored adaptive LASSO expectile method. For a Cox regression, \cite{Li.Pak.Todem.20} performs the automatic variable selection via a penalized nonparametric maximum likelihood estimator with an adaptive LASSO penalty. \cite{Li.Shi.Fan.Zhong.24} proves the oracle properties of the penalized weighted LS estimator with censored outcome adaptive LASSO penalty for survival outcomes which include AFT and Cox proportional hazards models. All these methods will be difficult to apply if the sample size is very large. In this case, we propose in the present paper, based on these estimators, aggregated estimators which will have the same theoretical properties as the estimator calculated on the full database. \\
  The paper is structured as follows. Section \ref{section_model} introduces the right-censored model, the general assumptions and notations. The building of the aggregated censored adaptive LASSO estimators are also presented. In Section \ref{section_estimators}, four of these estimators are defined and studied asymptotically, while in Section \ref{section_choixBIC} we present BIC-type criteria for selecting the tuning parameter which is used in the expression of the aggregated estimator.  In Section \ref{section_simus} the results of a numerical study are carried out. The proofs of the theoretical results are relegated in Section \ref{section_proofs}.
\section{Model and assumptions}
\label{section_model}
In this section, after introducing the accelerated failure time model, we present general assumptions and notations used throughout the paper. Finally, the construction of the aggregated censored adaptive LASSO estimator is presented as a useful estimator when the number of observations is very large and when automatic selection of  relevant variables is also required for a large number of explanatory variables. Let us start with some notations which will be used throughout the paper. All vectors are column; matrices and vectors are denoted with a boldface. For a set ${\cal A}$, we denote by $|{\cal A}|$ or $Card \{\cal A\}$ its cardinality, by ${\cal A}^c$ the complementary set. We shall use the notation $[x]$ for the integer part of the real number $x$.\\
Consider the following linear model on  $n$ observations:
\begin{equation}
\label{eq1}
T^*_i=\eX_i^\top \eb+\varepsilon_i, \qquad i=1, \cdots , n, 
\end{equation}
where $\eX_i$ is a  random vector of $p$ observable explanatory variables,  $\varepsilon_i$ is a  continuous unobservable random variable of errors,  $\eb \in \R^p$ is the vector of parameters and $\ebo=(\beta^0_1, \cdots , \beta^0_p)$ its true  (unknown) value.  The random variable  $T_i=\exp(T_i^*)$, $i=1, \cdots , n$, is the failure time (or survival time).\\
Let   ${\cal C}_i$ be the censoring variable (censoring time) for the $i$th observation.  In this paper we assume that $T_i$ is randomly right-censored by ${\cal C}_i$, such that the censoring time $T_i$ cannot always be observed. Thus, due to censoring, the observed variables are $(Y_i,\eX_i,\delta_i)_{1 \leqslant i \leqslant n}$, with $Y_i \equiv \min(T_i,{\cal C}_i)$ and the failure indicator $\delta_i \equiv \e1_{T_i \leq {\cal C}_i}$ which indicates whether $T_i$ was observed or not. We thus have an accelerated failure time (AFT) model.\\
We denote by $\PP_\eX$ the probability law of the random vector $\eX$ and by $\eE_\eX$ the expectation with respect to the distribution of $\eX$.  Similarly, we denote by $\PP_{\cal C}$, $\PP_\varepsilon$, the probability laws of $\cal C$ and $\varepsilon$, respectively and by $\PP$ the joint probability of $(\eX, {\cal C}, \varepsilon)$, and $\eE$ the corresponding expectation. \\
In the following, throughout this paper we will denote by $c$ a generic constant, without interest, which does not depend on $n$. We also denote by $Y$, ${\cal C}$, $T$, $\eX$, the generic variable for $Y_i$, ${\cal C}_i$, $T_i$, $\eX_i$. The components of  $\eX$ are $(X_1, \cdots , X_p)$ and those of  $\eX_i$ are $(X_{1i}, \cdots , X_{pi})$.\\
For the censoring variable ${\cal C}$, for all  $t>0$, we define its survival function   $G_0(t) \equiv \PP_{\cal C}[{\cal C} >t]$  and $1-G_0(t)$ its distribution function. \\
Since in applications, the survival function $G_0$  of the censoring variable $\cal C$ is often  unknown, its Kaplan-Meier estimator can be considered:
\[
\widehat G_n(t) \equiv \prod^n_{ \substack{
		i=1\\ Y_i \leq t}} \bigg( \frac{n -R_i}{n-R_i+1}\bigg)^{\delta_i},
\]
with $R_i$ the rank of $Y_i$ in $(Y_i)_{1 \leqslant i \leqslant n}$.\\

We now introduce the classical assumptions  on  ${\cal C}_i$, $\varepsilon_i$, $\eX_i$, $T_i$, for any $i=1, \cdots , n$, that will be considered for all the estimation methods presented in Section \ref{section_estimators}, where three further assumptions will be considered, taking into account the specificity of each method.  
\begin{description}
	\item \textbf{(A1)} The  censoring variable ${\cal C}_i$ is independent of  $\varepsilon_i$  and  $\eX_i$. The variable  ${\cal C}_i$ is also independent of the  failure time $T_i$  conditional on $\eX_i$.
\item \textbf{(A2)} The random vectors $(T_i,{\cal C}_i,\eX_i)_{ 1 \leqslant i \leqslant n}$ are independent and identically distributed (i.i.d.).
\item  \textbf{(A3)} The random variables $(\varepsilon_i)_{ 1 \leqslant i \leqslant n}$ are i.i.d.
\item  \textbf{(A4)} $\PP[t \leq T \leq {\cal C}] \geq \zeta_0 >0$ for all $t \in [0,B]$, with  $\zeta_0$ a positive  constant and  constant $B$  the maximum follow-up.
\end{description}
For the survival function $G_0$ we consider the classical condition:
\begin{description}
	\item \textbf{(A5)} The survival function $G_0$ is continuous and its derivative is uniformly bounded on $[0,B]$.
\end{description}
The design $(\eX_i)_{1 \leqslant i \leqslant n}$ and $\ebo$ satisfy the following two assumptions:
\begin{description}
\item \textbf{(A6)} The vector $\eX$  of explanatory variables is with compact support  and the true parameter $\ebo$ belongs to the interior of a known compact set.
	\item \textbf{(A7)} The matrix $\lim_{n \rightarrow \infty} n^{-1} \sum^n_{i=1} \eE [\eX_i \eX_i^\top/G_0(Y_i)]$ converges to a  defined positive matrix.
	%	\item \textbf{(C.3)}  $(\varepsilon_i)_{1 \leqslant i \leqslant n}$ sont indépendantes et la densité $f$ satisfait la condition $0 \leq \sup_{s}f(s| \eX=\ex) \leq B_0$.
\end{description}
Note that assumptions (A4) and (A5) ensure the consistency of the Kaplan-Meier estimator $\widehat G_n$.\\
\medskip
Moreover, we obviously have $\eE_{\cal C}[\delta_i] =G_0(Y_i)$ for any $i=1, \cdots ,n$.\\  
Assumptions (A1)-(A5) are commonly  considered in literature for the censored models (the reader is advised to consult \cite{Ciuperca.24}'s article for detailed references). Assumptions (A6), (A7) are also standard in the right-censored models (see for example \cite{De Backer.19}, \cite{Shows.Lu.Zhang.10}, \cite{Tang-Zhou-Wu.12}, \cite{Johnson.09}, \cite{Yu.Ye.Wang.24}). All these assumptions are necessary to prove the consistency and asymptotic normal properties of censored estimators.\\
We are going to work in the case when the number  $p$ of the explanatory variables  can be large, in which case we should perform automatic variable selection. For this purpose, the parameter $\eb$  is estimated by different adaptive LASSO methods that must be chosen according to the  assumptions made on $\varepsilon$.  \\
For the AFT model of relation (\ref{eq1}), let us introduce the index set of true non-zero coefficients: 
 $$
 {\cal A} \equiv \{j \in \{ 1, \cdots ,p\}; \; \beta^0_j \neq 0\}.
 $$
 Since $\ebo$ is unknown, so is the set ${\cal A}$. Without reducing the generality, we assume that  ${\cal A}$ contains  the first  $q \equiv |{\cal A}|$ natural numbers: ${\cal A}=\{1, \cdots , q\}$.  So its complementary set is ${\cal A}^c=\{q+1, \cdots , p\}$ and then $\ebo=(\eb^0_{\cal A},\oo_{|{\cal A}^c|})$.\\
  
 In order to study the properties of the estimators proposed in the following section, for $z, t \in [0,B]$, $j=1, \cdots , n$, consider the following random processes:
\begin{equation*}
%	\left\{  
	\begin{split}
		y(z)  &\equiv \lim_{n \rightarrow \infty} \frac{1}{n} \sum^n_{i=1} \e1_{Y_i \geq z}, \\
		M_j^{\cal C}(t) & \equiv (1- \delta_j) \e1_{Y_j \leq t} - \int^t_0 \e1_{Y_j \geq z} d \Lambda_{\cal C}(z),\\
		\Lambda_{\cal C}(t)& \equiv -\log(G_0(t)),
	\end{split}
%	\right. 
\end{equation*}
where $\Lambda_{\cal C}$ is the cumulative hazard function of the censoring  variable ${\cal C}$. \\
%Let us remark that $\{ M_j^{\cal C}(t)\}$ is a  martingale with respect to the $\sigma$-filtration: $\sigma \{\e1_{Y_j \geq u}, (1-\delta_j) \e1_{Y_j \leq u}, 0 \leq u \leq t; \eX_j; j=1, \cdots, n\}$ (see \cite{Sun.Zhang.09}).\\
Let us consider  $f_\varepsilon$ the density of $\varepsilon$ and for $\tau \in (0,1)$, $b_\tau \in \R$ the quantile of order $\tau$ of $\varepsilon$: $\PP_\varepsilon[\varepsilon \leq b_\tau]=\tau$. The  true value $b^0_\tau$ of $b_\tau$ is supposed  unknown.\\
We denote by  $\eX_{{\cal A}}$ the sub-vector of $\eX$ with indices in  ${\cal A}$ and for the $i$-th observation  $\eX_{{\cal A},i} \equiv (X_{ji})_{j \in {\cal A}}$. For a vector $\eb$, we use the notational conventions $\eb_{\cal A}$ for its sub-sector containing the corresponding components of ${\cal A}$.  
Consider for $b \in \R$ the following notations 
\begin{equation}
	\label{lee}
	\left\{  
	\begin{split}
		\eS_{{\cal A},b}& \equiv \eE_\eX[f_\varepsilon(b|\eX)\eX_{{\cal A}} \eX^\top_{{\cal A}}],\\
		\eS_{{\cal A}} &\equiv \eS_{{\cal A},0}.
			\end{split}
	\right. 
\end{equation}
In this paper, the number $n$ of observations can be very large, leading to numerical processing difficulties. To overcome the challenge of processing  all the data while at the same time automatically selecting the significant variables, we will share the observations $1, \cdots, n$ in $K$ groups, with $K=o(n)$, each group containing $n/K$ observations. Without reducing generality, we assume that $n$ is a multiple of $K$ and that each group has the same number of observations, $N \equiv n/K$. Thus, in the $k$-th group, with $k \in \{1, \cdots , K\}$, we consider the following observations ${\cal U}_k \equiv \{k, K+k, 2K+k, \cdots, (N-1)K+k\}$. This judicious  sharing of observations in the  $K$ groups allows us to keep the evolution of the survival curve $G_0(t)$ in each group $k$.\\
The same assumption $K=o(n)$ was considered in the paper of \cite{Su.Yin.Zhang.Zhao.23}. 
If the number  $n$ of observations is very large, corroborated with the possibility that $p$ is too, calculate the censored adaptive LASSO estimators may be impossible, due to the limited memory capacity of the computer. This is why the observations will be divided into several blocks on which estimators will be computed and which will then be  used to calculate an aggregated estimator, which will have the same asymptotic properties as the initial estimator (but which is impossible to compute). \\
In the next section, on observations  $1, \cdots , n$ we first calculate a censored estimator, denoted by $\widetilde{\eb}_n$ which will afterwards  be used to calculate a censored adaptive LASSO censored estimator, denoted  $\widehat{\eb}_n$. Numerical calculation difficulties  of these estimators can also  arise when $n$ is very large. Thus, for  group $k$ of observations ${\cal U}_k$ we calculate the same estimators, noted $\widetilde{\eb}^{(k)}_n=\big(\widetilde \beta^{(k)}_{n,1}, \cdots , \widetilde \beta^{(k)}_{n,p}\big)$  and $\widehat{\eb}^{(k)}_n=\big(\widehat \beta^{(k)}_{n,1}, \cdots , \widehat \beta^{(k)}_{n,p}\big)$, respectively.   For the estimators $\widehat{\eb}_n=\big(\widehat\beta_{n,1}, \cdots , \widehat\beta_{n,p}\big)$ and $\widehat{\eb}^{(k)}_n=\big(\widehat\beta_{n,1}^{(k)}, \cdots , \widehat\beta_{n,p}^{(k)}\big)$ consider the following  sets:
\[\widehat {\cal A}_n \equiv \{j \in \{1, \cdots, p\}; \; \widehat \beta_{n,j} \neq 0\}\]
and 
\[
\widehat {\cal A}^{(k)}_n \equiv \{j \in \{1, \cdots , p\}; \; \widehat\beta_{n,j}^{(k)} \neq 0\}.
\]
These two sets are estimators of ${\cal A}$.\\
 In order to study the oracle properties of $\overset{\vee}{\eb}_{n}$, we consider that the sequence of tuning parameters  $(\lambda_n)_{n \in \N}$ is such that 
\begin{equation}
	\label{eln}
	\lambda_n \rightarrow \infty, \quad \lambda_n=o(n^{1/2}).
\end{equation}
The tuning parameter $\lambda_{{\cal U}_k}$ depends on the observations of ${\cal U}_k$, for $k=1, \cdots ,K$. \\
For a fixed natural number   $w \in \{1, \cdots , K-1\}$, based on these estimators and sets, we will consider the following set of indices:
\[
\overset{\vee}{\cal A}_n \equiv \{j \in \{1, \cdots, p\}; \; \sum^K_{k=1}\e1_{\widehat \beta^{(k)}_{n,j}}> w\},
\]
which will be used to build the censored aggregated censored adaptive LASSO  estimator: 
\begin{equation}
	\label{agr}
\overset{\vee}{\eb}_n \equiv \frac{\bfQ_{ \overset{\vee}{\cal A}_n}}{K} \sum^K_{k=1} \widehat\eb_{n,\overset{\vee}{\cal A}_n}^{(k)},
\end{equation}
with the $p$-square matrix:
\[
\bfQ \equiv \textrm{diag} \bigg( \e1_{\sum^K_{k=1} \e1_{\widehat{\beta}^{(k)}_{n,1} \neq 0}\geq w  }, \e1_{\sum^K_{k=1} \e1_{\widehat{\beta}^{(k)}_{n,2} \neq 0}\geq w  }, \cdots , \e1_{\sum^K_{k=1} \e1_{\widehat{\beta}^{(k)}_{n,p} \neq 0}\geq w  }\bigg).
\]
%\[
%\bfQ \equiv \textrm{diag}(e_1, \cdots, e_p)
%, \qquad \textrm{where } \qquad  e_j \equiv \left\{
%\begin{array}{lll}
%	1, & & \textrm{if } \displaystyle{ \sum^K_{k=1} \e1_{\widehat{\beta}^{(k)}_{n,j} \neq 0}\geq w  },\\
%	0, & & \textrm{otherwise}.
%\end{array}
%\right.
%\] 
Based on  the matrix $\bfQ$ we construct   $  \textbf{Q}_{ \overset{\vee}{\cal A}_n }$ which is a  matrix   of dimension $p \times | \overset{\vee}{\cal A}_n |$, containing the columns of $\bfQ$ corresponding to the indices of  $ \overset{\vee}{\cal A}_n$. \\
This way of composing $\overset{\vee}{\cal A}_n$ and $\overset{\vee}{\eb}_n $ yields estimators that satisfy the oracle properties, i.e. sparsity and asymptotic normality for the estimators of non-zero coefficients. These estimators will be defined in the next section, considering four estimation techniques. \\
In this paper we consider that the numbers $p$ of explanatory variables and $|{\cal A}|$ of non-zero coefficients, do not depend on $n$.\\
Since $K=o(n)$, together with the  Central Limit Theorem (CLT)  for independent Bernoulli random variables (see \cite{Ciuperca.21}), we deduct that $w=O(K^{1/2})$.  
Then, from now on we will assume $w=O(K^{1/2})=o(n^{1/2})$.
\section{Estimators }
\label{section_estimators}
For each of the four estimation methods presented in this section, we first define the censored adaptive LASSO  estimators $\widehat \eb_n$ calculated on all $n$ observations. After presenting the asymptotic properties of $\widehat \eb_n$, we study the corresponding aggregate estimator $\overset{\vee}{\eb}_n$ calculated by (\ref{agr}), with the estimators $\widehat\eb_{n,\overset{\vee}{\cal A}_n}^{(k)}$ corresponding to following techniques: median, quantile, expectile and LS.  The main results for the four aggregated censored adaptive LASSO estimators are presented in next four subsections, while the proofs are all postponed to Section \ref{section_proofs}.
 \subsection{Median method}
 \label{meth_LAD}
 If $F_\varepsilon$ is the distribution function of $\varepsilon$, then in this subsection we assume that  $F_\varepsilon(0)=1/2$ and the estimation method is the median one. 
So, in this subsection we will study censored estimators obtained by the median method.\\
Let be the following random $|{\cal A}|$-vectors:
\begin{equation*}
%	\left\{  
	\begin{split} 
		\eh_{\cal A}(z)& \equiv \lim_{n \rightarrow \infty} \frac{1}{n} \sum^n_{i=1} \frac{\delta_i}{G_0(Y_i)} \e1_{Y_i \geq z} \eX_{{\cal A},i} \big(\e1_{\varepsilon_i <0} - \e1_{\varepsilon_i >0}\big), \qquad z \in [0,B],\\
		\es_i  &\equiv \frac{\delta_i}{G_0(Y_i)} \eX_{{\cal A},i} \big(\e1_{\varepsilon_i <0} - \e1_{\varepsilon_i >0}\big)+\int^B_0 \frac{\eh_{\cal A}(z)}{y(z)} dM^{\cal C}_i(z), \qquad i=1, \cdots, n.
		%,\\
		%\etha_{\cal A} & \equiv \bigg(\frac{sign(\beta^0_1)}{\beta^0_1}, \cdots,\frac{sign(\beta^0_{|{\cal A}|})}{\beta^0_{|{\cal A|}}} \bigg).
	\end{split}
%	\right. 
\end{equation*}
The censored adaptive LASSO median  estimator $\widehat \eb_n$ of $\ebo$ defined by \cite{Shows.Lu.Zhang.10} is the  minimizer of the following process:
 \begin{equation}
 	\label{bn_aLAD}
\sum^n_{i=1}\frac{\delta_i}{\widehat G_n(Y_i)}  \left| \log(Y_i) - \eX_i^\top \eb\right| +\lambda_n \sum^p_{j=1} \frac{|\beta_j|}{|\widetilde \beta_{n,j}|},
 \end{equation}
with $\widetilde{\eb}_n=(\widetilde{\beta}_{n,1}, \cdots ,\widetilde{\beta}_{n,p} )$ the censored median   estimator:
\begin{equation}
	\label{bn_LAD}
\widetilde{\eb}_n \equiv \argmin_{\eb \in \R^p} \sum^n_{i=1}\frac{\delta_i}{\widehat G_n(Y_i)}  \left| \log(Y_i) - \eX_i^\top \eb\right| .
\end{equation}
 and $(\lambda_n)_{n \in \N}$ a set of tuning parameters. In order to study the asymptotic properties of  the estimators considered in this subsection, let us formulate an additional assumption on the distribution of $\varepsilon$, in addition to those given in Section \ref{section_model}.
\begin{description}
	\item \textbf{(A8)}  The density $f_\varepsilon$ of model error  $\varepsilon$ satisfies the properties that there are two positive constants $b_0$ and $B_0$ such that $f_\varepsilon(0|\eX=\ex) \geq b_0 >0$, $|f'_\varepsilon((0|\eX=\ex)| \leq B_0$, $\sup_{s}f_\varepsilon(s|\eX=\ex) \leq B_0$.
\end{description}
Assumption (A8) was considered by \cite{Shows.Lu.Zhang.10} for showing the consistency and oracle properties of $\widehat \eb_n$. 
Under assumptions (A1)-(A8), by the proof of  Theorem 2 (i) and (ii) of \cite{Shows.Lu.Zhang.10} we have that, if $\lambda_n$ satisfies (\ref{eln}), then $\widehat \eb_n $ can be considered, with probability converging to 1, as $ \left( \ebo_{\cal A}+n^{-1/2}\eu_{\cal A}, \oo_{|{\cal A}^c|}\right)$, with $\eu_{\cal A} \in \R^{|{\cal A}|}$, $\|\eu_{\cal A} \|_2 \leq c < \infty$. By the same proof we have that the $|{\cal A}|$-vector $\eu_{\cal A}$ is the minimizer of
\begin{equation*}
	\frac{1}{2} \eu^\top_{\cal A} \eS_{\cal A} \eu_{\cal A}+\frac{\eu^\top_{\cal A}}{\sqrt{n}} \sum^n_{i=1}\es_i .
\end{equation*}
Then
\begin{equation}
	\label{eq2}
	\eu_{\cal A}=-\eS_{\cal A}^{-1} \frac{1}{\sqrt{n}}\sum^n_{i=1} \es_i ,
\end{equation}
which implies
\begin{equation}
	\label{eq3}
	\big(\widehat \eb_n - \ebo\big)_{\cal A} =- n^{-1/2}\eS_{\cal A}^{-1} \frac{1}{\sqrt{n}}\sum^n_{i=1} \es_i(1+o_\PP(n^{-1/2})).
\end{equation}
On the other hand, taking into account assumptions (A1)-(A4), we have by the central limit theorem (CLT): 
\begin{equation}
	\label{nb}
n^{-1/2} \sum^n_{i=1} \es_i  \overset{\cal L} {\underset{n \rightarrow \infty}{\longrightarrow}} {\cal N}_{|{\cal A}|}(\oo_{|{\cal A}|}, \eE[\es_1 \es_1^\top]).
\end{equation}
Note that under assumptions (A1)-(A8), by Theorem 2(i) of \cite{Shows.Lu.Zhang.10}, we have $\PP [\widehat{\cal A}_n ={\cal A}] {\underset{n \rightarrow \infty}{\longrightarrow}} 1$.\\
On observations $ {\cal U}_k$ (that is $k$th group of observations) corresponding to  relations (\ref{bn_aLAD}) and (\ref{bn_LAD}), the following censored estimators  can be calculated, respectively:
\begin{equation*}
%	\left\{  
	\begin{split} 
		\widehat \eb_n^{(k)}& \equiv \argmin_{\eb \in \R^p} \bigg( \sum_{i \in {\cal U}_k} \frac{\delta_i}{\widehat G_n(Y_i)} \big|\log (Y_i) - \eX_i^\top \eb\big|+\lambda_{{\cal U}_k}\sum^p_{j=1} \frac{|\beta_j|}{ |\widetilde \beta^{(k)}_{n,j}|}\bigg),\\
		\widetilde \eb_n^{(k)}& \equiv \argmin_{\eb \in \R^p}  \sum_{i \in {\cal U}_k} \frac{\delta_i}{\widehat G_n(Y_i)} \big|\log (Y_i) - \eX_i^\top \eb\big|.
	\end{split}
%	\right. 
\end{equation*}
Note also that the index set $\widehat {\cal A}^{(k)}_n $ satisfies the property: $\lim_{n \rightarrow \infty}\PP [\widehat {\cal A}^{(k)}_n=\widehat{\cal A}_n ={\cal A}] = 1$. \\
By the following theorem we show that the aggregated censored adaptive LASSO median estimator $\overset{\vee}{\eb}_{n}$  calculated by relation (\ref{agr}) satisfies the sparsity property as the unaggregated estimator $\widehat \eb_n$, the minimizer of  (\ref{bn_aLAD}). This means that the both estimators $\overset{\vee}{\eb}_{n}$ and $\widehat \eb_n$ select, with probability converging towards 1, the indices of the non-zero coefficients, i.e. of the relevant variables. Moreover, we show that the two estimators $ \overset{\vee}{\eb}_{n,{\cal A}}$ and $ \widehat{\eb}_{n,{\cal A}}$ of the non-zero coefficients have an identical asymptotic normal distribution. \\
Consequently, on a massive database, calculating  $\overset{\vee}{\eb}_{n}$ on $K$ subgroups of observations allows us to obtain an estimator with the same oracle proprieties as the estimator that would be calculated on the full database. On the other hand, the convergence rate of $\overset{\vee}{\eb}_{n}$ towards $\ebo$ depends on the number $K$ of groups, more precisely it is of order $(K/n)^{1/2}$. So the convergence rate may be slower that that of $\widehat \eb_n$, which is of order $n^{-1/2}$.
\begin{theorem}
	\label{th_oracle}
	Under assumptions (A1)-(A8), if moreover   $K=o(n)$, $w=O(K^{1/2})$, $F_\varepsilon(0)=1/2$, $(\lambda_n)_{n \in \N}$ satisfies (\ref{eln}), then:\\
	(i) $\lim_{n \rightarrow \infty} \PP\big[ \overset{\vee}{\cal A}_n= {\cal A} \big]=1$,  $\lim_{n \rightarrow\infty} \PP[ \overset{\vee}{\eb}_{n, {\cal A}^c}=\textbf{0}_{|{\cal A}^c|}]=1$ and  $\|  \overset{\vee}{\eb}_{n,{\cal A}} -\eb^0_{\cal A}\|_2=O_{\PP} \big( (  {K}/{n})^{1/2} \big)$.\\
	(ii) $n^{1/2} \big(  \overset{\vee}{\eb}_{n,{\cal A}} -\eb^0_{\cal A} \big) \overset{\cal L} {\underset{n \rightarrow \infty}{\longrightarrow}}  {\cal N} \pth{ \textbf{0}_{|{\cal A}|}, \eS^{-1}_{\cal A} \eV_{\cal A} \eS^{-1}_{\cal A}}$, with the matrices $\eV_{\cal A} \equiv \eE[\es_1 \es_1^\top]$ and $\eS_{\cal A}$ defined in relation (\ref{lee}).
\end{theorem}
%The proof of Theorem \ref{th_oracle} is given in Section \ref{section_proofs}.
\subsection{Quantile method}
\label{meth_quantile}
In this subsection, we will generalize the median method presented in subsection \ref{meth_LAD}. \\
For a fixed quantile index  $\tau \in (0,1)$, consider the check function $\rho_\tau:\R \rightarrow \R_+$ defined, for $u \in \R$ by, 
$$\rho_\tau(u)  = u [\tau - \e1_{u \leq 0}].$$
Then, in this subsection, for model errors  $\varepsilon$ we suppose that $F_\varepsilon(0)=\PP_\varepsilon[\varepsilon <0] =\tau$, with $\tau $ known.\\
Let us first consider the censored adaptive LASSO quantile estimator  $ \widehat \eb_n$ introduced and studied by \cite{Tang-Zhou-Wu.12} as the  minimizer of the following process:
\begin{equation}
	\label{bn_aquant}
	  \sum^n_{i=1}  \frac{\delta_i}{\widehat G_n(Y_i)} \rho_{\tau} (Y_i-\eX_i^\top \eb) +\lambda_n \sum^p_{j=1} \frac{|\beta_j|}{|\widetilde \beta_{n,j}|},
\end{equation}
with $\widetilde{\eb}_n=\big(\widetilde{\beta}_{n,1}, \cdots, \widetilde{\beta}_{n,p}\big)$ the censored quantile estimator defined by:
\begin{equation}
	\label{bn_quant}
   \widetilde{\eb}_n   \equiv \argmin_{ \eb }  \sum^n_{i=1}  \frac{\delta_i}{\widehat G_n(Y_i)} \rho_{\tau} (Y_i -\eX_i^\top \eb).
\end{equation}
In order to present the asymptotic properties of $\widehat \eb_n$ and those of its aggregated estimator, let us consider the deterministic value:
\[
	\varphi_0 \equiv \eE_\eX[f_\varepsilon(0|\eX)], 
\]
and the random vectors:
\begin{equation*}
%	\left\{  
	\begin{split} 
		%	 \eD_j & \equiv \eE_\eX[f(b^0_{\tau_j}|\eX) \eX \eX^\top]\\
	%	\eD_{\cal A}& \equiv   \eS_{{\cal A},0} \\
		\es_{i}^{(1)}  &\equiv    \frac{\delta_i}{G_0(Y_i)} \eX_{{\cal A},i} \big(\e1_{\varepsilon_i <0} - \tau\big),\\
		\es_{i}^{(2)}  &\equiv    \frac{\delta_i}{G_0(Y_i)} \frac{\eX_{{\cal A},i}}{n} \big(\e1_{\varepsilon_i <0} - \tau\big)\sum^n_{k=1} \int^B_0 \frac{\e1_{Y_i \geq z}}{y(z)} dM^{\cal C}_k(z).
	\end{split}
%	\right. 
\end{equation*}
Given the proof of Theorem 2.3(ii) of \cite{Tang-Zhou-Wu.12}, under assumptions (A1)-(A8), (\ref{eln}),  then  the estimator  $\widehat{\eb}_n$ can be written, with probability converging to 1, as $\ebo+n^{-1/2} \eu$, with $\eu=(\eu_{\cal A}, \oo_{|{\cal A}^c|})$, $\eu_{\cal A} \in \R^{|{\cal A}|}$, $\| \eu_{\cal A}\|_2 \leq c< \infty$.\\
Hence, with the matrix  $\eS_{{\cal A}}$ defined in relation  (\ref{lee}), the $|{\cal A}|$-vector $\eu_{\cal A}$ is the minimizer of the following process (see the proof of Theorem 2.3 of \cite{Tang-Zhou-Wu.12}), with probability converging to 1:
\[
\frac{1}{2}\eu_{\cal A}^\top \eS_{{\cal A}} \eu_{\cal A}+\eu_{\cal A}^\top \frac{1}{ \sqrt{n}}\sum^n_{i=1} \es_{i}^{(1)}+\eu_{\cal A}^\top \frac{1}{ \sqrt{n}}\sum^n_{i=1} \es_{i}^{(2)} ,
\] 
from where
\begin{equation}
	\label{eq4}
	\eu_{\cal A}=- \big(\eS_{{\cal A}}\big)^{-1} n^{-1/2}\sum^n_{i=1}\big(\es_{i}^{(1)}+\es_{i}^{(2)}\big),
\end{equation}
with 
\begin{equation}
	\label{eq_etoile}
	\left\{  
	\begin{split} 
	 n^{-1/2} \sum^n_{i=1} \es_{i}^{(1)} \overset{\cal L} {\underset{n \rightarrow \infty}{\longrightarrow}} \eW_1 \sim  {\cal N}_{|{\cal A}|}(\oo_{|{\cal A}|},...),\\
	n^{-1/2} \sum^n_{i=1} \es_{i}^{(2)} \overset{\cal L} {\underset{n \rightarrow \infty}{\longrightarrow}} \eW_2 \sim  {\cal N}_{|{\cal A}|}(\oo_{|{\cal A}|},...).
	\end{split}
	\right. 
\end{equation} 
Then, for the censored adaptive LASSO quantile estimator, the relation similar to (\ref{eq3}) is 
\begin{equation}
	\label{r123}
	\big(\widehat \eb_n - \ebo\big)_{\cal A} =- n^{-1/2}\eS_{{\cal A}}^{-1} \bigg(n^{-1/2}\sum^n_{i=1} \big(\es_{i}^{(1)}+\es_{i}^{(1)}\big) \bigg)(1+o_\PP(n^{-1/2})),
\end{equation}
On the other hand, on observations ${\cal U}_k$ ($k$h group of observations) the censored adaptive LASSO quantile estimators (corresponding to  relations (\ref{bn_aLAD}) and (\ref{bn_LAD})) are respectively:
\begin{equation*}
%	\left\{  
	\begin{split} 
		\widehat \eb_n^{(k)}& \equiv \argmin_{\eb \in \R^p}    \sum_{i \in {\cal U}_k} \frac{\delta_i}{\widehat G_n(Y_i)} \rho_{\tau}\big(\log (Y_i)- \eX_i^\top \eb\big)+\lambda_{{\cal U}_k}\sum^p_{j=1} \frac{|\beta_j|}{ |\widetilde \beta^{(k)}_{n,j}|}, \\
		\widetilde \eb_n^{(k)}& \equiv \argmin_{\eb \in \R^p}   \sum_{i \in {\cal U}_k} \frac{\delta_i}{\widehat G_n(Y_i)}\rho_{\tau}\big(\log (Y_i) - \eX_i^\top \eb\big).
	\end{split}
%	\right. 
\end{equation*}
With these elements, we can show  by the following theorem that the aggregated censored adaptive LASSO quantile estimator $\overset{\vee}{\eb}_{n}$ calculated  by relation (\ref{agr}) satisfies the oracle properties.
\begin{theorem}
	\label{th_oracleq}
Under assumptions (A1)-(A8), if moreover   $K=o(n)$, $w=O(K^{1/2})$, $F_\varepsilon(0)=\tau$, $(\lambda_n)_{n \in \N}$ satisfies (\ref{eln}), then:\\
	(i) $\lim_{n \rightarrow \infty} \PP\big[ \overset{\vee}{\cal A}_n= {\cal A} \big]=1$,  $\lim_{n \rightarrow\infty} \PP[ \overset{\vee}{\eb}_{n, {\cal A}^c}=\textbf{0}_{|{\cal A}^c|}]=1$ and  $\|  \overset{\vee}{\eb}_{n,{\cal A}} -\eb^0_{\cal A}\|_2=O_{\PP} \big( (  {K}/{n})^{1/2} \big)$.\\
	(ii) $n^{1/2} \big(  \overset{\vee}{\eb}_{n,{\cal A}} -\eb^0_{\cal A} \big) \overset{\cal L} {\underset{n \rightarrow \infty}{\longrightarrow}}  {\cal N} \pth{ \textbf{0}_{|{\cal A}|}, \eS_{{\cal A}}^{-1} \Var[\eW_1+\eW_2] \eS_{{\cal A}}^{-1} }$.
\end{theorem}
%The proof of Theorem \ref{th_oracleq} is given in Section \ref{section_proofs}.\\

We will further generalize these results by considering, as in \cite{Tang-Zhou-Wu.12}, the  quantile index  $\tau_j=j/(1+J)$ for $j\in \{1, \cdots , J\}$  and $b^0_{\tau_1}, \cdots , b^0_{\tau_J}$ the true (unknown) values of the corresponding quantiles. The natural number $J$ is known and does not depend on $n$. \\
In this case, the censored adaptive LASSO composite quantile estimator  $\big(\widehat b_{n,1}, \cdots , \widehat b_{n,J}, \widehat \eb_n\big)$ minimizes in respect to $\big(b_1, \cdots , b_J, \eb \big)$ the following random process (see \cite{Tang-Zhou-Wu.12}):
\begin{equation}
	\label{bn_aquantJ}
\sum^J_{j=1} \sum^n_{i=1}  \frac{\delta_i}{\widehat G_n(Y_i)} \rho_{\tau_j} (Y_i-b_j-\eX_i^\top \eb) +\lambda_n \sum^p_{j=1} \frac{|\beta_j|}{|\widetilde \beta_{n,j}|},
\end{equation}
with $\widetilde{\eb}_n=\big(\widetilde{\beta}_{n,1}, \cdots, \widetilde{\beta}_{n,p}\big)$ the censored quantile estimator defined by:
\begin{equation}
	\label{bn_quantJ}
\big(\widetilde b_{n,1}, \cdots, \widetilde b_{n,J},\widetilde{\eb}_n \big) \equiv \argmin_{b_1, \cdots, b_J,\eb } \sum^J_{j=1} \sum^n_{i=1}  \frac{\delta_i}{\widehat G_n(Y_i)} \rho_{\tau_j} (Y_i-b_j-\eX_i^\top \eb).
\end{equation}
Let be the deterministic value: $ \varphi_j \equiv \eE_\eX[f_\varepsilon(b^0_{\tau_j}|\eX)]$, $ j=1, \cdots , J$, the deterministic matrix:
\[
 \eD_{\cal A} \equiv \sum^J_{j=1} \eS_{{\cal A},b^0_{\tau_j}}
\]
and the random $|{\cal A}|$-vectors:
\begin{equation*}
%	\left\{  
	\begin{split} 
			\es_{i,J}^{(1)}  &\equiv \sum^J_{j=1}  \frac{\delta_i}{G_0(Y_i)} \eX_{{\cal A},i} \big(\e1_{\varepsilon_i <b^0_{\tau_j}} - \tau_j\big),\\
	\es_{i,J}^{(2)}  &\equiv \sum^J_{j=1}  \frac{\delta_i}{G_0(Y_i)} \frac{\eX_{{\cal A},i}}{n} \big(\e1_{\varepsilon_i <b^0_{\tau_j}} - \tau_j\big)\sum^n_{k=1} \int^B_0 \frac{\e1_{Y_i \geq z}}{y(z)} dM^{\cal C}_k(z),
	\end{split}
%	\right. 
\end{equation*}
with $\eS_{{\cal A},b^0_{\tau_j}}$ defined in relation (\ref{lee}).
Then, the following theorem, which is a generalization of Theorem \ref{th_oracleq}, holds.
\begin{theorem}
	\label{th_oracleqJ}
Under assumptions (A1)-(A8),   if moreover   $K=o(n)$, $w=O(K^{1/2})$,  $(\lambda_n)_{n \in \N}$ satisfies (\ref{eln}), then:\\
	(i) $\lim_{n \rightarrow \infty} \PP\big[ \overset{\vee}{\cal A}_n= {\cal A} \big]=1$,  $\lim_{n \rightarrow\infty} \PP[ \overset{\vee}{\eb}_{n, {\cal A}^c}=\textbf{0}_{|{\cal A}^c|}]=1$ and  $\|  \overset{\vee}{\eb}_{n,{\cal A}} -\eb^0_{\cal A}\|_2=O_{\PP} \big( (  {K}/{n})^{1/2} \big)$.\\
	(ii) $n^{1/2} \big(  \overset{\vee}{\eb}_{n,{\cal A}} -\eb^0_{\cal A} \big) \overset{\cal L} {\underset{n \rightarrow \infty}{\longrightarrow}}  {\cal N} \pth{ \textbf{0}_{|{\cal A}|}, \eD^{-1}_{\cal A} \Var[\eW_{1,J}+\eW_{2,J}] \eD^{-1}_{\cal A}}$, with the normal random vectors $\eW_{1,J}$ and $\eW_{2,J}$ defined in the proof given in Section \ref{section_proofs}.
\end{theorem}

Obviously, if  $J=1$ and $\tau_1=0.5$ we obtain the median method for which, instead of assumption $F_\varepsilon(0)=0.5$, we take $F_\varepsilon(b_1)=0.5$, with $b_1$ unknown and to be estimated.
 \subsection{Expectile method}
 \label{meth_expectile}
 Since the quantile function is not derivable, in this subsection, we will consider the expectile estimation method. 
 For  the expectile index $\tau \in (0,1)$, the expectile function is defined by:
 \[
 \rho_\tau(u)=|\tau - \e1_{u <0}| u^2, \qquad u \in \R,
 \]  
 which is derivable. The non-derivability of the quantile function can poses certain numerical problems and also difficulties in the theoretical study of the estimators. As will be seen from the simulations in some cases the quantile estimators will perform less well than expectile ones. \\
 The derivative of $\rho_\tau(u)$   is $g_\tau(u) \equiv\rho'_\tau(u-t)|_{t=0}= - 2 \tau u \e1_{u \geq 0}- 2(1-\tau)u \e1_{u<0}$ and the second derivative is $h_\tau(u) \equiv \rho''_\tau(u-t)|_{t=0}=2 \tau \e1_{u \geq 0}+2(1-\tau) \e1_{u<0}$. \\
 
For the model error $\varepsilon$ we consider the following assumption:
 \begin{description}
 	\item \textbf{(A9)}   $\eE_\varepsilon[\varepsilon^4] < \infty$ and $\eE_\varepsilon[g_\tau(\varepsilon)]=0$
 \end{description}
 and for the design $(\eX_i)_{1 \leqslant i \leqslant n}$:
  \begin{description}
 \item  \textbf{(A10)} $n^{-1} \sum^n_{i=1} \eX_i \eX_i^\top \overset{\PP_\eX} {\underset{n \rightarrow \infty}{\longrightarrow}}  \eE_\eX [\eX \eX^\top]$, with $ \eE_\eX [\eX \eX^\top]$ a  positive definite matrix.
  \end{description}
Then, in this subsection we consider assumptions (A1)-(A7), (A9) and (A10). Note that if $\tau=1/2$ then $\eE_\varepsilon[g_\tau(\varepsilon)]=0$ in assumption (A9) becomes $\eE_\varepsilon[\varepsilon]=0$, which is the standard condition considered for the LS estimation method. Note that assumptions (A9) and (A10) are classic for the expectile method (see \cite{Ciuperca.24}, \cite{Gu-Zou.16}).\\
 
 For $z \in [0,B]$, let be the following random $p$-vector 
 $$\ek(z) \equiv \lim_{n \rightarrow \infty} \frac{1}{n} \sum^n_{i=1}  \frac{\delta_i}{G_0(Y_i)}  \e1_{Y_i \geq z}g_\tau(\varepsilon_i) \eX_i ,$$
 which, by assumption (A4), is bounded on  $[0,B]$,  with the constant $B$ defined in assumption  (A4).\\
 We also define the random $|{\cal A}|$-vector:
 \begin{equation}
 	\label{kaz}
  \ek_{\cal A}(z) \equiv \lim_{n \rightarrow \infty} \frac{1}{n} \sum^n_{i=1}  \frac{\delta_i}{G_0(Y_i)}  \e1_{Y_i \geq z}g_\tau(\varepsilon_i) \eX_{{\cal A},i} .
  \end{equation}
The censored expectile estimator was defined and studied by \cite{Ciuperca.24}:
 \begin{equation}
 	\label{bn_expectile}
 \widetilde \eb_n \equiv \argmin_{\eb \in \R^p} \sum^n_{i=1} \frac{\delta_i}{\widehat G_n(Y_i)} \rho_\tau(\log(Y_i)- \eX_i^\top \eb ),
 \end{equation}
 the components of the $p$-vector $\widetilde \eb_n $ being denoted $\big(\widetilde \beta_{n,1}, \cdots, \widetilde \beta_{n,p}\big)$.   
\cite{Ciuperca.24} still defines and studies the asymptotic properties of the censored adaptive LASSO expectile estimator:
\begin{equation}
	\label{ebn}
\widehat \eb_n \equiv \argmin_{\eb \in \R^p} \bigg(  \sum^n_{i=1}\frac{\delta_i}{\widehat G_n(Y_i)}\rho_\tau(\log(Y_i)-\eX_i^\top \eb)+ \lambda_n \sum^p_{j=1} \widehat \omega_{n,j}|\beta_j| \bigg),
\end{equation}
with the adaptive weights $ \widehat \omega_{n,j}\equiv |\widetilde \beta_{n,j}|^{-\gamma} $, $\widetilde \beta_{n,j}$ the $j$-th component of $\widetilde \eb_n$ given by (\ref{bn_expectile})  and $\gamma >0$ a known parameter.  \\
Let be the random $|{\cal A}|$-vector 
$$\evs_{\cal A} \equiv \frac{\delta}{  G_0(Y)} g_\tau(\varepsilon)  \eX_{{\cal A}}+\int^B_0 \frac{\ek_{\cal A}(z)}{y(z)}dM^{\cal C}(z) 
$$
 and the $|{\cal A}|$-square matrix 
 $$
 \eE[\evs_{\cal A} \evs^\top_{\cal A}]=\textbf{S}_{1,{\cal A}}+\textbf{S}_{2,{\cal A}},
 $$ 
with matrix $\textbf{S}_1\equiv \eE_\varepsilon\big[ g^2_\tau(\varepsilon)\big] \eE_\eX\big[ \eX \eX^\top/ G_0(Y)\big]$ and $\textbf{S}_2 \equiv \eE \big[\int^B_0 {\ek(z) \ek^\top(z)}/{y(z)}d \Lambda_C(z) \big]$. Let also be  $p$-square matrix  $\textbf{S}_3 \equiv  \eE_\varepsilon \big[h_\tau(\varepsilon)\big] \eE_\eX [\eX \eX^\top]$. \\
Before stating the following result, let us consider the random $|{\cal A}|$-vector:
 \begin{equation}
 	\label{sit}
 	\widetilde \es_i \equiv \frac{\delta_i}{G_0(Y_i)}g_\tau(\varepsilon_i) \eX_{{\cal A},i}+\int^B_0 \frac{\ek_{\cal A}(z) }{y(z)}dM^{\cal C}_i(z).
 \end{equation}
 With these notations we can prove that the aggregated censored adaptive LASSO expectile estimator $ \overset{\vee}{\eb}_{n} $ of (\ref{agr}), with $\widehat\eb_{n,\overset{\vee}{\cal A}_n}^{(k)}$ calculated by  (\ref{ebn}), satisfies the oracle properties. 
 \begin{theorem}
 	\label{th_oraclee}
 	Under assumptions (A1)-(A7), (A9) and (A10), if moreover   $K=o(n)$, $w=O(K^{1/2})$, $F_\varepsilon(0)=1/2$, $(\lambda_n)_{n \in \N}$ satisfies (\ref{eln}), then:\\
 	(i) $\lim_{n \rightarrow \infty} \PP\big[ \overset{\vee}{\cal A}_n= {\cal A} \big]=1$,  $\lim_{n \rightarrow\infty} \PP[ \overset{\vee}{\eb}_{n, {\cal A}^c}=\textbf{0}_{|{\cal A}^c|}]=1$ and  $\|  \overset{\vee}{\eb}_{n,{\cal A}} -\eb^0_{\cal A}\|_2=O_{\PP} \big( (  {K}/{n})^{1/2} \big)$.\\
 	(ii) $n^{1/2} \big(  \overset{\vee}{\eb}_{n,{\cal A}} -\eb^0_{\cal A} \big) \overset{\cal L} {\underset{n \rightarrow \infty}{\longrightarrow}}  {\cal N} \pth{ \textbf{0}_{|{\cal A}|}, \eS_{{\cal A}}^{-1} \Var[\eW_1+\eW_2] \eS_{{\cal A}}^{-1} }$, where the  $|{\cal A}|$-square matrix $\eSS$ is defined by  $ \eSS \equiv \eE^{-2}_\varepsilon[h_\tau(\varepsilon)]\eE_\eX[\eX_{{\cal A}}\eX^\top_{{\cal A}}]^{-1}   \big(\textbf{S}_{1,{\cal A}}+\textbf{S}_{2,{\cal A}}\big) \eE_\eX[\eX_{{\cal A}}\eX^\top_{{\cal A}}]^{-1}$.
 	% and $\eo^0_{\cal A} \equiv \lim_{n \rightarrow \infty} \widehat \eo_{n,{\cal A}} =|\eb^0_{\cal A}|^{-\gamma}$. 
 \end{theorem}
Note that we can write $\eSS=\textbf{S}^{-1}_{3,{\cal A}}(\textbf{S}_{1,{\cal A}}+\textbf{S}_{2,{\cal A}})\textbf{S}^{-1}_{3,{\cal A}}$.  

\subsection{Least squares method}
\label{meth_LS}
If in  subsection \ref{meth_expectile} we consider the particular case for the expectile index $\tau=1/2$, then, $\rho_{1/2}(x)=x^2/2$, $g_{1/2}(x)=-x$ and $h_{1/2}(x)=1$, from where  $\eE_\varepsilon[ g_{1/2}(\varepsilon)]=\eE_\varepsilon[  \varepsilon]=0$, $\eE_\varepsilon\big[ g^2_{1/2}(\varepsilon)\big]=\Var[\varepsilon]$, $\eE_\varepsilon[h_{1/2}(\varepsilon)] =1$.  In this case, for $z\in [0,B]$, the random vector $\ek_{\cal A}(z)$ of relation (\ref{kaz}) becomes:
$$
\ek_{\cal A}(z) \equiv - \lim_{n \rightarrow \infty} \frac{1}{n} \sum^n_{i=1}  \frac{\delta_i}{G_0(Y_i)}  \e1_{Y_i \geq z} \varepsilon_i \eX_{{\cal A},i}  ,
$$
while $\widetilde \es_i $ of relation (\ref{sit}) becomes:
\[
\widetilde \es_i \equiv - \frac{\delta_i}{G_0(Y_i)} \varepsilon_i \eX_{{\cal A},i}+\int^B_0 \frac{\ek_{\cal A}(z) }{y(z)}dM^{\cal C}_i(z).
\]
Thus, taking into account relation (\ref{beA}) in the proof of Theorem \ref{th_oraclee}(\textit{ii}) for the censored adaptive LASSO expectile estimator, we obtain for the censored adaptive LASSO LS estimator $\widehat \eb_n$ the following relation:
\[
(\widehat \eb_n -\ebo)_{\cal A} =-n^{-1/2}\eE_\eX[\eX_{\cal A} \eX_{{\cal A}}^\top] ^{-1} n^{-1/2} \sum^n_{i=1}\widetilde \es_i \big(1+ o_\PP(n^{-1/2})\big).
\]
The matrix  $\eSS$ of Theorem \ref{th_oraclee} becomes in this particular case:
\[
 \eE_\eX[\eX_{\cal A} \eX_{{\cal A}}^\top] ^{-1} \bigg(\eE_\varepsilon[\varepsilon^2] \eE_\eX \bigg[\frac{\eX_{\cal A} \eX_{{\cal A}}^\top}{G_0(Y)}\bigg]+ \eE \bigg[\int^B_0 \frac{{\ek_{\cal A}(z) \ek_{\cal A}^\top(z)}}{{y(z)}}d \Lambda_C(s) \bigg]\bigg) \eE_\eX[\eX_{\cal A} \eX_{{\cal A}}^\top] ^{-1}.
\]
For $\gamma =1$, the censored adaptive LASSO LS estimator of  (\ref{ebn}) becomes the estimator studied by \cite{Wang.Song.2011}.
\section{Choice of $\lambda_n$ by the BIC criterion}
\label{section_choixBIC}
In this section we present BIC-type criteria for selecting the tuning parameter $\lambda_n$. For each of the four methods presented in Section \ref{section_estimators}, we denote by $\widehat \eb_n(\lambda_n)$ the censored adaptive LASSO estimator obtained for a fixed $\lambda_n$.\\
For the censored adaptive LASSO median estimator $\widehat \eb_n(\lambda_n)$ minimizer of (\ref{bn_aLAD}) in  subsection \ref{meth_LAD}, \cite{Shows.Lu.Zhang.10} proposes as a choice for the tuning parameter $\lambda_n$ the one that minimizes the following BIC-type criterion:
\begin{equation}
	\label{BIC_LAD}
BIC(\lambda_n)=\frac{1}{BIC(0)} \sum^n_{i=1} \frac{\delta_i}{\widehat G_n(Y_i)} \big|\log(Y_i)-\eX_i^\top \widehat \eb_n(\lambda_n)\big|  +|\widehat{\cal A}_n(\lambda_n)|  \frac{\log n}{n},
\end{equation}
with
\[
BIC(0) \equiv  \sum^n_{i=1} \frac{\delta_i}{\widehat G_n(Y_i)} \big|\log(Y_i)-\eX_i^\top \widetilde\eb_n\big|.
\]
The estimator  $\widetilde \eb_n$ is calculated by (\ref{bn_LAD}).\\
From this criterion, we will first propose another, generalizing it for the quantile method presented  in  subsection \ref{meth_quantile}:
\begin{equation}
	\label{BIC_quantile}
BIC(\lambda_n)=\frac{1}{BIC(0)} \sum^n_{i=1}   \frac{\delta_i}{\widehat G_n(Y_i)} \rho_{\tau}\big(\log(Y_i)-\eX_i^\top \widehat \eb_n(\lambda_n)\big)  +|\widehat{\cal A}_n(\lambda_n)|  \frac{\log n}{n},
\end{equation}
with
\[
BIC(0) \equiv  \sum^n_{i=1}  \frac{\delta_i}{\widehat G_n(Y_i)} \rho_{\tau} \big(\log(Y_i)-\eX_i^\top \widetilde \eb_{n}\big).
\]
The censored quantile estimator $\widetilde \eb_n$ is calculated by  (\ref{bn_quant}). Note that  criterion (\ref{BIC_quantile}) is different from that proposed in the paper of \cite{Tang-Zhou-Wu.12} which considered: 
\[
BIC(\lambda_n)=\log \bigg( \frac{1}{J} \sum^J_{j=1} \frac{1}{n} \sum^n_{i=1}\frac{\delta_i}{\widehat G_n(Y_i)} \big|Y_i-\widehat b_j(\lambda_n) -\eX_i^\top \widehat \eb_{n} (\lambda_n)\big|\bigg)+\big| \widehat {\cal A}_n(\lambda_n)\big| \frac{\log n}{n}.
\]
For the expectile method presented in subsection  \ref{meth_expectile}, the optimal $\lambda_n$ is the minimizer of the following criterion:
\begin{equation}
	\label{BIC_expectile}
BIC(\lambda_n)=\frac{1}{BIC(0)} \sum^n_{i=1}  \frac{\delta_i}{\widehat G_n(Y_i)} \rho_{\tau}\big(\log(Y_i -\eX_i^\top \widehat \eb_n(\lambda_n)\big)  +|\widehat{\cal A}_n(\lambda_n)| \frac{\log n}{n},
\end{equation}
with 
\[
BIC(0)\equiv  \sum^n_{i=1}   \frac{\delta_i}{\widehat G_n(Y_i)} \rho_{\tau} \big(\log(Y_i) -\eX_i^\top \widetilde \eb_{n}\big).
\]
The estimator $ \widehat \eb_n(\lambda_n)$ which depends on  $\lambda_n$ is the minimizer of relation (\ref{ebn}) and $\widetilde \eb_n$ is calculated by (\ref{bn_expectile}).\\
Obviously we can also calculate for each $k$ group of observations   ${\cal U}_k$ the best tuning parameter $\lambda_{{\cal U}_k}$.\\

Note that an aggregated estimator in a censoring model with massive data has been proposed in paper of \cite{Su.Yin.Zhang.Zhao.23} but where the loss function is weighted least squares, and therefore derivable. Derivability makes it easy to carry out a numerical study of the method proposed by    \cite{Su.Yin.Zhang.Zhao.23}. In the present paper (except for the censored quantile estimator) the loss function is not derivable. On the other hand, the criterion proposed in \cite{Su.Yin.Zhang.Zhao.23} is different from ours, i.e. that of relation (\ref{BIC_expectile}) for $\tau=1/2$. The only point in common with our BIC criterion is the penalty $|\widehat{\cal A}_n(\lambda_n)|  \frac{\log n}{n}$. 
\section{Simulation study}
\label{section_simus}
In this section we perform a numerical study by Monte Carlo simulations to evaluate the properties of the aggregated censored adaptive LASSO estimators $  \overset{\vee}{\eb}_{n} $  and compare them with the censored adaptive LASSO estimators $\widehat \eb_n$ calculated on the  full  data. The BIC criterion is also studied. All simulations are performed using R language. The used packages are {\it SALES, quantreg} and {\it cqrReg} for adaptive LASSO expectile, quantile, composite quantile methods, respectively. \\
We will examine an important property of the aggregated estimator $ \overset{\vee}{ \eb}_n$, namely its sparsity (or automatic selection). For this reason, we consider an estimation of a coefficient to be a false zero if the true value of the coefficient is different from 0 and the estimate is 0. On the other hand, for a true value of a zero coefficient estimated by a non-zero value, we say that we have a false non-zero. More precisely, for $M$ Monte Carlo replications we calculate:
\[
\%\textrm{ of false zeros}=\frac{100}{M} \sum^M_{m=1} \frac{Card \{j \in {\cal A}; \; \overset{\vee}{ \beta}^{(m)}_{n,j}=0\}}{ |{\cal A}|}
\]
and 
\[
\%\textrm{ of false non-zeros}= \frac{100}{M} \sum^M_{m=1} \frac{Card \{j \in {\cal A}^c; \; \overset{\vee}{ \beta}^{(m)}_{n,j}\neq 0\}}{ |{\cal A}^c|}, 
\]
with $\overset{\vee}{ \beta}^{(m)}_{n,j} $ the aggregated censored adaptive LASSO  estimation of $\beta_j$ obtained at the $m$ Monte Carlo replication, for $m=1, \cdots , M$.\\

Throughout this section, the following model is considered:
\begin{equation}
\label{eq1s}
T_i^*=\beta^0_0+\sum^p_{j=1}\beta^0_j X_{ji}+\varepsilon_i, \qquad i=1, \cdots, n,
\end{equation}
where, unless otherwise stated, the design is $X_{ji} \sim {\cal N}(1,1)$ for any $i=1, \cdots , n$ and $j=1, \cdots , p$. The censoring variable is with uniform distribution ${\cal C}_i \sim {\cal U}[0,c_1]$, with the constant $c_1$ chosen such that the censoring rate is $25\%$. The model errors are standard Gumbel $\varepsilon \sim {\cal G}(0,1)$.\\
In  subsections \ref{subsec_BIC}-\ref{subsec et veb} we consider the set ${\cal A}=\{1, 2\}$, with $\beta^0_1=1$, $\beta^0_2=-2$.   \\
For the expectile method, the considered power of the weights in the adaptive penalty is $\gamma=1$. On the other hand, the expectile index $\tau$ must satisfy condition $\eE[g_\tau (\varepsilon)]=0$ of assumption (A9) and then, for a given distribution of $\varepsilon$, we will consider an estimate for $\tau$. Thus, the empirical estimation considered for the expectile  index $\tau$ is:
\begin{equation}
	\label{ee}
\widehat \tau_n^{(E)}=\frac{n^{-1}\sum^n_{i=1} \varepsilon_i \e1_{\varepsilon_i <0}}{n^{-1} \big( \sum^n_{i=1} \varepsilon_i \e1_{\varepsilon_i <0}-\sum^n_{i=1} \varepsilon_i \e1_{\varepsilon_i >0}\big)}.
\end{equation}
  Similarly, an empirical estimation for quantile index $\tau$, such that $F_\varepsilon(0)=\tau$ is:
  \begin{equation*}
  \widehat \tau_n^{(Q)} =\frac{1}{n} \sum^n_{i=1} \varepsilon_i \e1_{\varepsilon_i <0}.
  \end{equation*}

\subsection{BIC criterion study}
\label{subsec_BIC}
\begin{figure}[h!] 
	\begin{tabular}{cc}
		\includegraphics[width=0.45\linewidth,height=4.5cm]{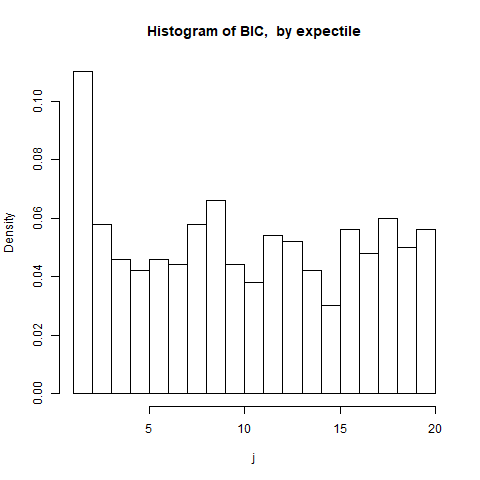} &
		\includegraphics[width=0.45\linewidth,height=4.5cm]{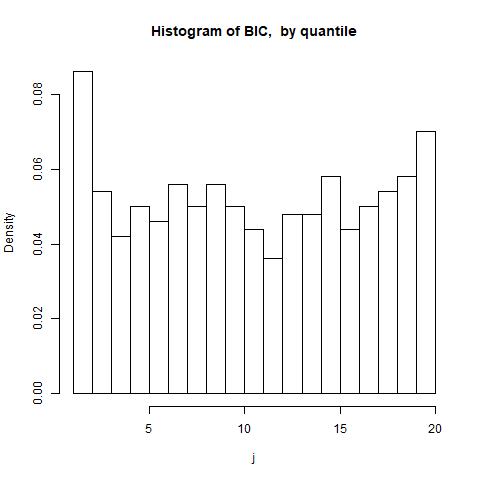} \\
		{\small (a)  $p=10$,   expectile method.} &
		{\small (b) $p=10$,  quantile method.}\\
		& \\
		\includegraphics[width=0.45\linewidth,height=4.5cm]{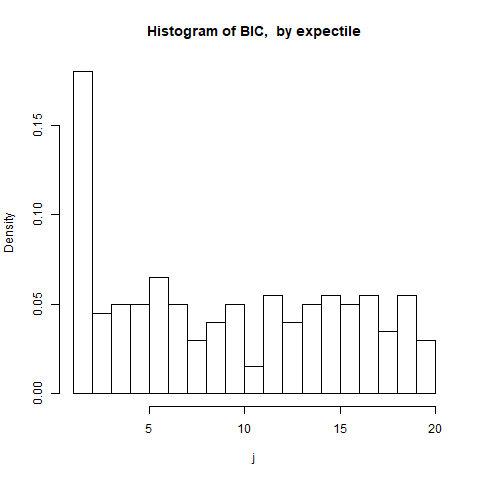} &
		\includegraphics[width=0.45\linewidth,height=4.5cm]{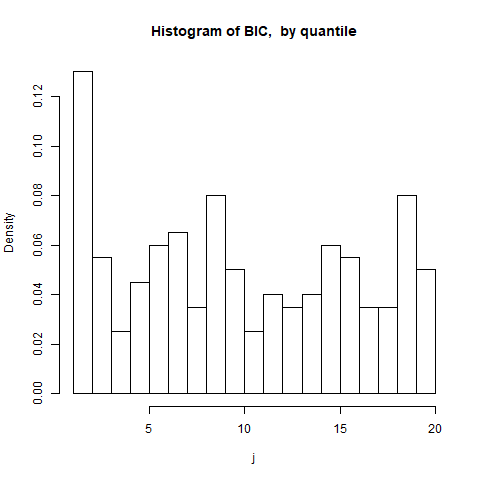} \\
		{\small (c) $p=100$,  expectile method.} &
		{\small (d) $p=100$,  quantile method.}
	\end{tabular}
	\caption{\small The histograms of indices  $j \in \{1, \cdots , 20\}$  which minimize the BIC criteria for a tuning parameter of the form   $\lambda_n=n^{1/2-1/(10j)}$, for $n=1000$, censored adaptive LASSO expectile and quantile methods.	}
	\label{fig_BIC_n2000}
\end{figure}

In this subsection we study the BIC criterion proposed in Section \ref{section_choixBIC} for choosing the tuning parameter $\lambda_n$. First of all, we consider two possibilities for the penalty in the BIC criterion. Specifically, let $\nu \equiv Card\{\delta_i=1; \; i=1, \cdots , n\}$ be, i.e. the number of times when the variable $T_i$ has been observed. In relations (\ref{BIC_LAD}), (\ref{BIC_quantile}), (\ref{BIC_expectile})  of BIC,  we will study the results obtained by choosing  $(\log n)/n$ or $(\log  \nu)/\nu$ in the penalty   (see Table \ref{Tab_BIC}).\\
We consider the following values for the sample size $n \in \{10^2, 10^3, 10^4\}$, the tuning parameter $\lambda_n=n^{1/2-1/(10j)}$ with $j \in \{ 1, 2, \cdots, 20\}$ and the number of parameters $p \in \{10, 50, 100, 500\}$. \\
We observe that the percentage of false zeros is always 0, while the percentage of false non-zeros does not vary with $n$ or $p$ but with the estimation method. The expectile method produces slightly fewer false non-zeros than the median method, which yields less than that quantile. There are also a few more false non-zeros when we take  $\nu$ in the BIC penalty. In view of these results, from now, on in all simulations we will take $(\log n)/n$ in the BIC penalty. \\
 \begin{table}
	\begin{center}
		\begin{tabular}{|B|B|B|BBB|BBB|}\hline  
			  $n$  	&    $p$   &    $n$ or $\nu$   & \multicolumn{3}{B|}{   $\%$ of false 0  }  & \multicolumn{3}{B|}{   $\%$ of false non-zeros}  \\ 
			\cline{4-9}
			& & &    expectile   &    median   &    quantile   &    expectile     &    median   &    quantile    \\ \hline 
			    100   &    10   &    $\nu$   &    0  &    0  &    0   &     0.02   &    1.11   &    1.93   \\
			& &    $n$   &    0  &    0  &    0  &    0.18   &    1  &    1.31  \\  \cline{2-9}
			&    50   &    $\nu$   &    0  &    0  &    0  &    2.75  &    9.68    &    10.49  \\
			& &    $n$   &    0  &    0  &    0   &    3.63   &    8.14   &    9.24   \\	   \hline
		   1000   &    10   &   $\nu$   &    0  &    0  &    0  &    0.05   &   0.05   &    1.65  \\
			 & &    $n$   &    0  &    0  &    0  &    0.02   &   0.02   &    1.47  \\  
		 	\cline{2-9}
		 	 &    100   &    $\nu$   &    0  &    0  &    0  &    0.07   &    0.13   &    0.98   \\ 
		 	 & &    $n$   &    0  &    0  &    0  &    0.27   &    0.06   &    0.71  \\ 
		 	 \cline{2-9}
		 	 &    500   &    $\nu$   &    0  &    0  &    0  &    0.36   &    0.80   &    1.53   \\ 
		 	 & &    $n$   &    0  &    0  &    0   &    0.65   &    0.45   &    1.23  \\  \hline
		 	    10000   &    10   & $\nu$ &    0  &    0  &    0  &    0   &    0.31   &   7.8   \\
		 	  & &    $n$   &    0  &    0  &    0  &    0.02   &    0.25   &   6.25  \\  \cline{2-9}
		 &    100   &    $\nu$   &    0  &    0  &    0   &    0  &    0.02   &   0.85   \\
		 & &    $n$   &    0  &    0  &    0  &    0   &    0.02   &   0.59   \\	   \hline
		\end{tabular} 
	\end{center}
	\caption{\small False zero and false non-zero detection rates for the corresponding BIC criteria presented in Section \ref{section_choixBIC}, when   $\beta^0_0=0$, ${\cal A}=\{1,2\}$, $\beta^0_1=1$, $\beta^0_2=-2$.}
	\label{Tab_BIC} 
\end{table}
In Figure \ref{fig_BIC_n2000} we present the histogram of $j \in \{1, \cdots , 20\}$,  which minimizes the BIC criterion for a tuning parameter of the form  $\lambda_n=n^{1/2-1/(10j)}$, for $n=1000$,  censored adaptive LASSO expectile and quantile methods  and $p \in \{10, 100\}$.
We obtain that the BIC criterion achieves its minimum mainly for $\lambda_n=n^{1/2-1/10}$, by both methods.
\subsection{Comparison of $\widehat \eb_n$ and $\overset{\vee}{\eb}_{n}$, with respect to $K$, $p$ and $w$ }
In Table \ref{Tab_comp_betas}, by 200   Monte Carlo replications,   the estimators $\widehat \eb_n$ and $\overset{\vee}{\eb}_{n}$ calculated by four censored estimation methods (expectile, median, quantile and composite expectile) are compared by the percentage presentation  of the false zeros and   of the false non-zeros when the dimension $p$ of $\eb$, the number $K$ of groups,  the value threshold $w$ are varied. For the censored adaptive LASSO composite quantile estimator of  (\ref{bn_aquantJ}) we take $J=10$.\\
We obtain that for the two estimators $\widehat \eb_n$ and $\overset{\vee}{\eb}_{n}$,   the percentage of false zeros  is  0  for all considered values of $p$, $K$, $w$, for the four methods. Then the percentage of false non-zeros is what distinguishes the four methods for each estimator $\widehat \eb_n$ and $\overset{\vee}{\eb}_{n}$.\\
 Concerning $\widehat \eb_n$, the median method produces more false non-zeros followed by the quantile method with the rate always below $5\%$. Emphasize that $|\widehat{\cal A}_n|$ remains constant equal to 0 for the censored adaptive LASSO composite quantile estimator. If $p \in \{50, 100\}$, then for the censored adaptive LASSO expectile, median, quantile, composite quantile estimators   calculated on all observations, we obtained $\|(\widehat{\eb}_{n} -\ebo)_{\cal A} \|_1=0.33$, $0.48$, $0.35$ and $2.99$, respectively. So, the bias of $\widehat \eb_{n}$, for a fixed method, does not vary with $p$.
On the other hand, the percentage of false non-zeros decreases as $p$ increases for expectile, median and quantile methods.\\
 For  $\overset{\vee}{\eb}_{n}$, the best results are obtained by expectile method and the worst by composite quantile.   Moreover, if $w=1$ and $K \in \{20, 100\}$, then for $\overset{\vee}{\eb}_{n}$ there are more  false non-zeros than for $\widehat \eb_n$ (except expectile when $K=20$). If $w=[\sqrt{K}]$, for  expectile, median, quantile techniques, then for $\widehat \eb_n$ that there are more false non-zeros  that for $\overset{\vee}{\eb}_{n}$. The aggregated censored adaptive LASSO composite quantile estimator  provides a lot of false zeros, which is why it will be abandoned in the following simulations.  These results are complemented by those in Table \ref{Tab_comp_betas_bis}. From Table \ref{Tab_comp_betas_bis} we deduce that for the aggregated censored adaptive LASSO composite quantile  estimator in order to make an efficient automatic selection (false non-zero rate below $5\%$) of null coefficients, then either the design must be centered, or the number $K$ of groups must be small and the value of the non-zero coefficients must be very far from zero.
Returning to Table \ref{Tab_comp_betas}, we deduce that the false non-zero rate does not vary with $p$   and increases when $K$ increases. Moreover, for $w=K^{1/2}$,  we make less than $1\%$ of false non-zeros, whatever the values of $K$ and $p$ for expectile, median and quantile methods.  %Moreover, $|\widehat{\cal A}_n|$ remains constant equal to 0 for all three methods, as already observed in Table \ref{Tab_BIC}, but also by censored adaptive LASSO composite quantile estimator.
\begin{table}
	\begin{center}
		\begin{tabular}{|B|BBBB|B|B|BBBB|}\hline  
			$p$	& \multicolumn{4}{B|}{ $\widehat \eb_n$}& $K$ & $w$   & \multicolumn{4}{B|}{$\overset{\vee}{\eb}_{n}$  }  \\ 
			\cline{2-5} 	\cline{8-11} 
		&  \multicolumn{4}{B|}{$\%$ of false non-zeros} 	& &   & \multicolumn{4}{B|}{$\%$ of false non-zeros} \\
		\cline{2-5} 	\cline{8-11}
		& expectile   & median & quantile & comp quant	&  &  & expectile   & median & quantile & comp quant   \\ \hline 
		50 &   & &  &  & 2 & 1  & 0 & 0& 0  & 25  \\
		\cline{6-11}
		  &   & &  &  & 20 & 1  & 0.02 & 10.4& 1.41  & 100  \\
			\cline{7-11}
		  	&  0.21 &4.83 &0.31 & 0 &  & $[\sqrt{K}]$=4  & 0 &0.04   & 0  & 100 \\
			\cline{6-11}
		&  & &   &   & 100 & 1  & 7.3 & 93 & 68 & 100 \\
		\cline{7-11}
	&  &  &   & 	&  & $[\sqrt{K}]$=10  & 0&0.5  & 0& 100 \\
		\cline{1-11}
		100 &  &  &   &  & 20 & 1  & 0 & 6.49 & 0.71 & 97 \\
		\cline{7-11}
		& 0.08 & 2.77& 2.65 &0&  & $[\sqrt{K}]$=4  & 0&0  & 0 & 58  \\
		\cline{6-11}
	&  &  &    &  	& 100 & 1 & 7.96 &88 & 62  & 100   \\
		\cline{7-11}
	&   &  &  &  	&  &$ [\sqrt{K}]$=10 & 0 & 0.10 & 0 & 100 \\
			   \hline
		\end{tabular} 
	\end{center}
	\caption{\small Comparison of   $\widehat \eb_n$ and $\overset{\vee}{\eb}_{n}$ by calculating the detection rates  of  false non-zero,  for $n= 10^5$, $\lambda_n=n^{1/2-1/10}$, $\beta^0_0=0$.}
	\label{Tab_comp_betas} 
\end{table}
\begin{table}
\begin{center}
	\begin{tabular}{|B|B|B|BBBB|}\hline  
 $\eb_{\cal A}$	&$\eX_i$ &	 $K$    & \multicolumn{4}{B|}{$\%$ of false non-zeros }  \\ 
		\cline{4-7}  
		    &  	&    & expectile   & median & quantile & comp quant   \\ \hline 
		(10,-20) & ${\cal N}(1,1) $& 2 & 3 & 22 & 7 & 5    \\ 
		& & 20 & 2 & 64 & 24 & 99  \\ 	\cline{2-7}
		 & ${\cal N}(0,1) $& 2 & 0 & 0.2 & 0 & 0    \\ 
		& & 20 & 0 & 11 & 5 & 32  \\ 	\cline{1-7} 
		(1,-2) & ${\cal N}(1,1) $& 2 & 0 & 0.2 & 0 & 100   \\ 
		& & 20 &0 & 0 & 0 & 100  \\ 	\cline{2-7}
		& ${\cal N}(0,1) $& 2 & 0 & 0 & 0 & 0    \\ 
		& & 20 & 0 & 0 & 0 & 0  \\    
		\hline
	\end{tabular} 
\end{center}
\caption{\small False non-zeros percentages  for $\overset{\vee}{\eb}_{n}$,  when $n= 10^5$, $\lambda_n=n^{1/2-1/10}$, $p=50$, ${\cal A}=\{1,2\}$, $\beta^0_0=0$, $w=[\sqrt{K}]$.}
\label{Tab_comp_betas_bis} 
\end{table}

\subsection{Study of  $\overset{\vee}{\eb}_{n}$ with respect to  $K$ and $w$}
\label{subsec et veb}
\begin{figure}[h!] 
	\begin{tabular}{cc}
		\includegraphics[width=0.45\linewidth,height=4.5cm]{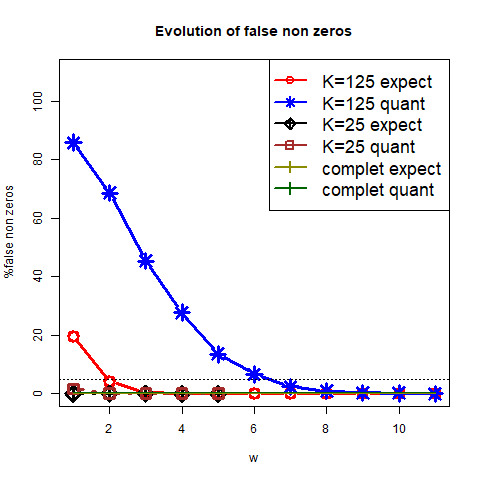} &
		\includegraphics[width=0.45\linewidth,height=4.5cm]{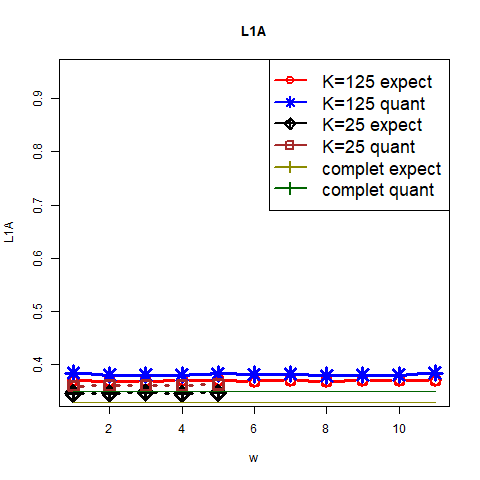} \\
		{\small (a) Evolution of false non-zeros.} &
		{\small (b) Evolution of $L1A=\|\big(\overset{\vee}{\eb}_{n}-\ebo \big)_{\cal A} \|_1$.}
	\end{tabular}
	\caption{\small Study of the estimator    $\overset{\vee}{\eb}_{n}$ with respect to values of $w$, for $K \in \{25, 125\}$, $p=50$. 	}
	\label{fig_evol_w}
\end{figure}

\begin{figure}[h!] 
	\begin{tabular}{cc}
		\includegraphics[width=0.45\linewidth,height=4.5cm]{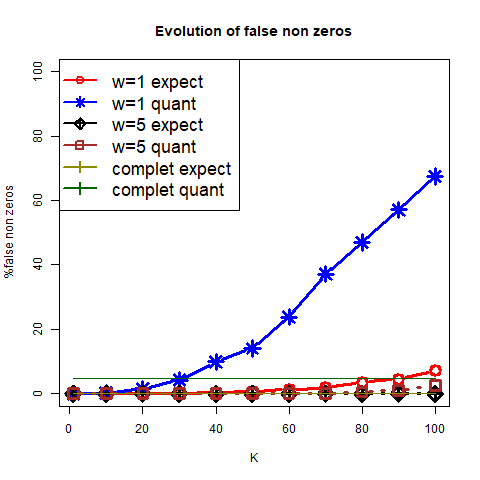} &
		\includegraphics[width=0.45\linewidth,height=4.5cm]{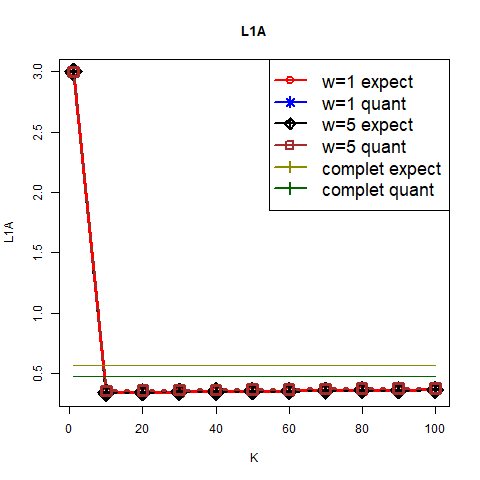} \\
		{\small (a) Evolution of false non-zeros.} &
		{\small (b) Evolution of $L1A=\|\big(\overset{\vee}{\eb}_{n}-\ebo \big)_{\cal A} \|_1$.}
	\end{tabular}
	\caption{\small Study of the estimator $\overset{\vee}{\eb}_{n}$   with respect to values of  $K$, for $w \in \{1, 5\}$,   $p=50$. }
	\label{fig_evol_K}
\end{figure}

\begin{figure}[h!] 
 	\begin{tabular}{cc}
 		\includegraphics[width=0.45\linewidth,height=4.5cm]{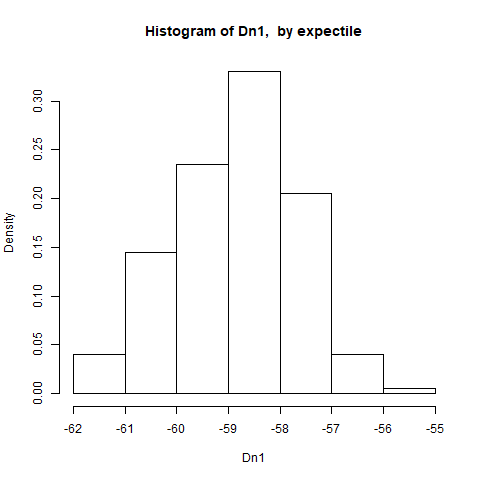} &
 		\includegraphics[width=0.45\linewidth,height=4.5cm]{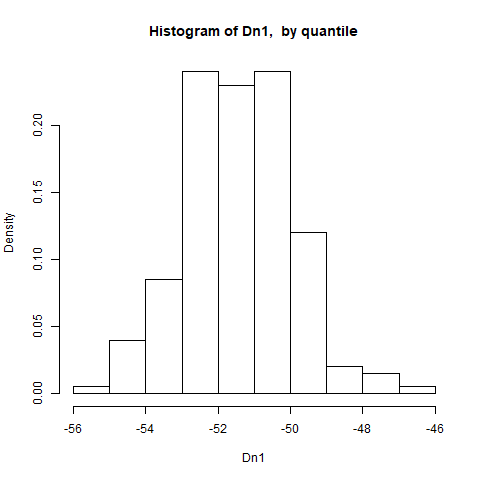} \\
 		{\small (a)  expectile method.} &
 		{\small (b)  quantile   method.}\\
 	\end{tabular}
 	\caption{\small The histogram of $Dn1 \equiv {\sqrt n}(\overset{\vee}{\beta}_{n,1}-\beta^0_1)$ by aggregated censored adaptive LASSO  methods.	}
 	\label{hist_normalite_beta1}
 \end{figure}

In this subsection we study the estimator $\overset{\vee}{\eb}_{n}$ with respect to  $K$ and $w$ for three methods: expectile, median and quantile. \\
From Table \ref{Tab_vbeta_K}, for 200 Monte Carlo replications of a model without intercept, $X_{ji} \sim {\cal N}(1,1)$ for any $i=1, \cdots , n$ and $j=1, \cdots , p$, we first deduce that for $p$ fixed ($p=50$),   the bias of $\overset{\vee}{\eb}_{n}$ and $\widehat{\eb}_{n}$ is about the same for a fixed estimation method and that by the median method the bias is greater. Biases are greater for a model with intercept ($\beta^0_0 \neq 0$) .  On the other hand, the median method produces more false non-zeros than the quantile, that does more than expectile, especially when $w=2$. If $w=[\sqrt{K}]$, then the percentage of false non-zeros is less than $1\%$ for a model without intercept and $100\%$ with intercept.\\
In Table \ref{Tab_vbeta_K_X0}, for a centered design, the bias of estimates  and false non-zero rates are smaller than those for a non-centered design.\\
In Table \ref{Tab_temps} we present execution times for a Monte Carlo replication. The simulations are performed on a computer Inter(R), Core (CM), 1.6 Ghz. Note that if $K=1$, then the estimations are based on full data.  The “Total” part includes data generation, data formatting, calculation of the Kaplan-Meier estimator $\widehat G_n$, calculation of $\overset{\vee}{\eb}_{n}$, $\overset{\vee}{\cal A}_n$ by the three methods.  For "expectile", "median",  "quantile"  we give only the computation time for $\overset{\vee}{\eb}_{n}$ by the expectile, median and quantile methods, respectively. We note that the execution time on full data is greater than the time on $K$ groups, by the three estimation methods. Furthermore, the execution time (for an entire Monte Carlo replication or for each of the estimation methods) decreases as $K$ increases up to a certain value and then starts to rise again with $K$.\\
In Figures \ref{fig_evol_w}(a) and (b), for the values of $K \in \{25, 125\}$, we plot the evolution of false non-zeros and of  $\|\big(\overset{\vee}{\eb}_{n}-\ebo \big)_{\cal A} \|_1$ as a function of $w$ values when $w \in \{1, 2, \cdots, [\sqrt{K}]\}$. In order to better situate the results presented in Figures \ref{fig_evol_w} and \ref{fig_evol_K}, let us recall that for $\widehat \eb_n$, the percentage of false non-zeros is respectively, 0.18, 4.98 and 0.56, $\|(\widehat{\eb}_{n} -\ebo)_{\cal A} \|_1=0.33$, $0.48$ and $0.35$ by the expectile, median and quantile methods, respectively. From Figure \ref{fig_evol_w}(a) we deduce that when $K=25$, then the $5\%$ rate of false non-zeros is reached for $w \geq 1$ for both estimation methods, while if $K$ is large (=125) then only for $w \geq 6$ do we have a rate lower than $5\%$ of false non-zeros for quantile method and $w \geq 2$ for expectile method.  We deduce that the more classes we take, the more we need to increase the value of $w$.   Hence, in order to expect the same rate of false non-zeros as on full data, for the quantile method the value of $w$ must be greater than by the expectile method. From Figure \ref{fig_evol_w}(b) we deduce that the bias of the aggregated estimator $\overset{\vee}{\eb}_{n}$ is approximately the same than that obtained on the full data by $\widehat \eb_n$ and it does not depend on the number of groups $K$ nor on the value of $w$ (results also corroborated by those in Table \ref{Tab_vbeta_K}).  \\
In Figures \ref{fig_evol_K}(a) and (b), for the values of $w \in \{1, 5\}$ we plot the evolution of false non-zeros and $\|\big(\overset{\vee}{\eb}_{n}-\ebo \big)_{\cal A} \|_1$ as a function of the number $K$ of groups. Figure \ref{fig_evol_K}(a) confirms that for quantile method a larger value for $w$ is required to avoid more than $5\%$ false non-zeros as the number $K$ of groups increases which is not the case for the expectile method, which gives a rate of less than $5\%$ for $K \leq 90$ and $w \geq 1$.  Moreover, for  given $K$ and $w$, we observe that  the expectile method produces fewer false non-zeros. For $K \leq 100$ and  $w=5$, then the aggregated expectile estimator makes as many false non-zeros as the estimator on the  full  data, while by the aggregated quantile estimator $\overset{\vee}{\eb}_{n}$ it is necessary to take for $w=5$ a value of $K \leq 80$ so that the false non-zeros are in the same proportion as those obtained for $\widehat \eb_n$. From Figure \ref{fig_evol_K}(b) we deduce that the bias increases very slightly as $K$ increases (in accordance with Theorem \ref{th_oracle}), but it does not depend on $w$, as already observed in Figure \ref{fig_evol_w}(a), and is very close to the bias of $\widehat \eb_n$.\\
For $K$ and $w$ fixed, we exemplify in Figures \ref{hist_normalite_beta1}(a) and (b) that Theorem \ref{th_oracle} is true for ${\sqrt n}(\overset{\vee}{\beta}_{n,1}-\beta^0_1)$ by the aggregated censored adaptive LASSO, expectile and quantile methods, respectively. Similar histograms were obtained (not presented for lack of space) for ${\sqrt n}(\overset{\vee}{\beta}_{n,2}-\beta^0_2)$.     The standard-deviations of ${\sqrt n}(\overset{\vee}{\beta}_{n,1}-\beta^0_1)$ obtained by 200 Monte Carlo replications are  1.16, and 1.50, by the expectile and quantile methods, respectively. Shapiro's test for normality of $\overset{\vee}{\beta}_{n,1}$ gives p-values 0.51 and 0.58 for the expectile and quantile methods, respectively.  Moreover, the standard-deviations of ${\sqrt n}(\overset{\vee}{\beta}_{n,2}-\beta^0_2)$ are 1, and 1.39, by the expectile and quantile methods, respectively. Shapiro's test for normality of $\overset{\vee}{\beta}_{n,2}$ gives p-values 0.14 and  0.87 for the expectile and quantile methods, respectively. These results are obtained for ${\cal A}=\{1,2\}$, $\beta^0_1=1$, $\beta^0_2=-2$, $K=100$, $n=10^5$, $p=50$ and $w=10$. 

\begin{table}
	\begin{center}
		\begin{tabular}{|B|B|B|BBB|BBB|}\hline  
$\beta^0_0$&	$w$ &	$K$ &  \multicolumn{3}{B|}{$\%$ of false non-zeros} &\multicolumn{3}{B||}{$\|(\overset{\vee}{\eb}_{n} -\ebo)_{\cal A} \|_1$}   \\ 
	\cline{3-9}  
	&	& & expectile   & median &quantile & expectile & median &quantile 	 \\ \hline 
0	&	2 &5  & 0  & 0.02 & 0& 0.33 & 0.48& 0.35 	 \\
			\cline{3-9}  
&&	25 &  0  & 3.9 &0.13 & 0.34 & 0.49& 0.36 \\ 	\cline{3-9}  
&&	50  & 0.03 & 24& 3.5 &	 0.35 & 0.50& 0.36 \\ \cline{3-9}  
&&	100	  & 1.07 & 81& 41&  0.36 & 0.51 & 0.37\\ \cline{2-9}  
&$[\sqrt{K}]$ &5  & 0   & 0.02 & 0& 0.33 & 0.48 &  0.35 \\ 	\cline{3-9}  
&&	25 &  0  & 0 &0 & 0.34 & 0.49& 0.36 \\ 	\cline{3-9}
&&  50 & 0   & 0.02 & 0 & 0.35 & 0.50 & 0.36 \\ 	\cline{3-9}  
&& 100 & 0   & 0.55 &0 & 0.36 & 0.51 & 0.37\\ 	\cline{1-9}  
2	&	2 &5  & 0  & 0.70 & 1& 0.70 & 1.27& 1.13 	 \\ \cline{3-9}
&&	50  & 100 & 100& 100 &	 0.38 & 2.06& 1.87 \\ \cline{2-9}    
	&	$[\sqrt{K}]$ &5  & 100  & 100 & 100& 0.92 & 2.27& 2.18 	 \\ \cline{3-9}
&&	50  & 100 & 100& 100 &	 0.38 & 2.06& 1.88
. \\ \cline{1-9}    
		\end{tabular} 
	\end{center}
	\caption{\small Study of $\overset{\vee}{\eb}_{n}$   when $n= 10^5$, $\lambda_n=n^{1/2-1/10}$,   $p=50$,  $\eX_i \sim {\cal N}(1,1)$.  Note that by 200   Monte Carlo,   when $\beta^0_0=0$, we obtain $\|(\widehat{\eb}_{n} -\ebo)_{\cal A} \|_1=0.33$, $0.48$ and $0.35$,  the percentage of false zeros is 0.18, 4.98 and 0.56, using the expectile, median and quantile methods, respectively.}
	\label{Tab_vbeta_K} 
\end{table}
\begin{table}
	\begin{center}
		\begin{tabular}{|B|B|B|BBB|BBB|}\hline  
			$\beta^0_0$&	$w$ &	$K$ &  \multicolumn{3}{B|}{$\%$ of false non-zeros} &\multicolumn{3}{B||}{$\|(\overset{\vee}{\eb}_{n} -\ebo)_{\cal A} \|_1$}   \\ 
			\cline{3-9}  
			&	& & expectile   & median &quantile & expectile & median &quantile 	 \\ \hline 
			0	&	2 &5  & 0  & 0 & 0& 0.35 & 0.39& 0.35 	 \\
			\cline{3-9}  
			&&	50  & 0 & 2.3& 0.29 &	 0.36 & 0.40& 0.37 \\  \cline{2-9}  
			&$[\sqrt{K}]$ &5  & 0   & 0 & 0& 0.35 & 0.39 &  0.36 \\  	\cline{3-9}
			&&  50 & 0   & 0 & 0 & 0.36 & 0.40 & 0.37 \\  	\cline{1-9}  
			2	&	2 &5  & 0  & 0 & 0& 0.13 & 0.84& 0.68 	 \\ \cline{3-9}
			&&	50  & 0.02 & 30& 4.64 &	 0.11 & 0.83& 0.67 \\ \cline{2-9}    
			&	$[\sqrt{K}]$ &5  & 0  & 0 & 0& 0.13 & 0.85& 0.69 	 \\ \cline{3-9}
			&&	50  & 0 & 0.06& 0 &	 0.11 & 0.83&  0.67\\ \cline{1-9}    
		\end{tabular} 
	\end{center}
	\caption{\small Study of $\overset{\vee}{\eb}_{n}$  when $n= 10^5$, $\lambda_n=n^{1/2-1/10}$,    $p=50$, $\eX_i \sim {\cal N}(0,1)$.  By 200   Monte Carlo replications, when $\beta^0_0=0$, we have $\|(\widehat{\eb}_{n} -\ebo)_{\cal A} \|_1=0.34$, $0.38$ and $0.35$, using the expectile, median and quantile methods and the percentage of false zeros is  0 for the three methods, respectively.}
	\label{Tab_vbeta_K_X0} 
\end{table}
\begin{table}
	\begin{center}
		\begin{tabular}{|B|B|BBBBBB|}\hline  
		Part of code executed & $p$ & $K=1$ & $K=5$ & $K=25$ & $K=50$ & $K=100$ & $K=200$ \\ \hline
		\textit{Total} & \textit{50} & \textit{31.9} & \textit{21.34} & \textit{ 13.92 }& \textit{10.29} & \textit{10.82}  & \textit{12.22} \\
		expectile & & 6.79 & 7.24 & 4.20 & 3.32  & 4.06 & 5.11 \\
		median & & 10.77 & 6.42 & 4.14 & 2.84 & 2.70 & 2.74 \\  
		quantile & & 13.28 & 6.61 & 4.18 & 2.99 & 2.75 & 2.99 \\ \hline
		\textit{Total }& \textit{100} & \textit{129.07} & \textit{69.66} & \textit{ 32.15} & \textit{29.48} & \textit{32.84} & \textit{29.09} \\
		expectile & & 19.38 & 15.56 & 8.13 & 9.40  & 12.37 & 13.86 \\
		median & & 59.97 & 24.72 & 11.23 & 9.21 & 9.27 & 6.60 \\ 
		quantile & & 47.85 & 27.75 & 11.15 & 9.28 & 9.59 & 6.97 \\\hline
		\end{tabular} 
	\end{center}
	\caption{\small Execution time (in seconds) for a Monte Carlo replication to calculate $\overset{\vee}{\eb}_{n}$ using the expectile, median and quantile methods, when $n= 10^5$, $\lambda_n=n^{1/2-1/10}$,   $p\in \{50, 100\}$, $w=5$. }
	\label{Tab_temps} 
\end{table}
 
\subsection{Automatic selection with respect to $\ebo$ values}
\label{subsec_ebo}
\begin{figure}[h!] 
	\begin{tabular}{cc}
		\includegraphics[width=0.45\linewidth,height=4.5cm]{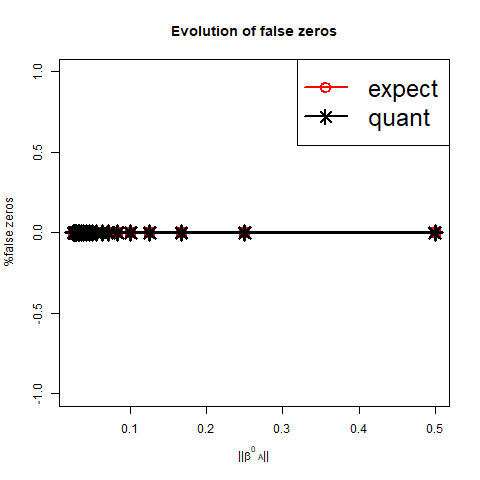} &
		\includegraphics[width=0.45\linewidth,height=4.5cm]{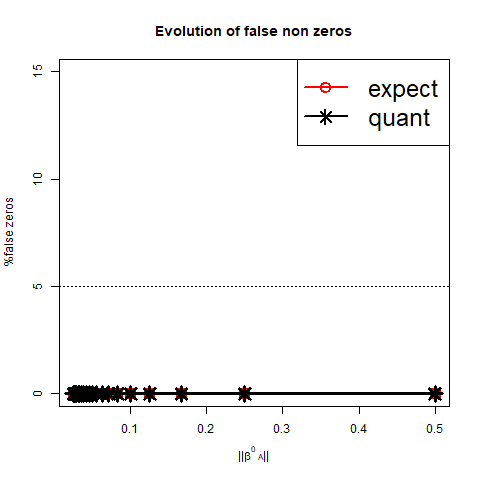} \\
		{\small (a) Evolution of false zeros.} &
		{\small (b) Evolution of false non-zeros.}\\
		\includegraphics[width=0.45\linewidth,height=4.5cm]{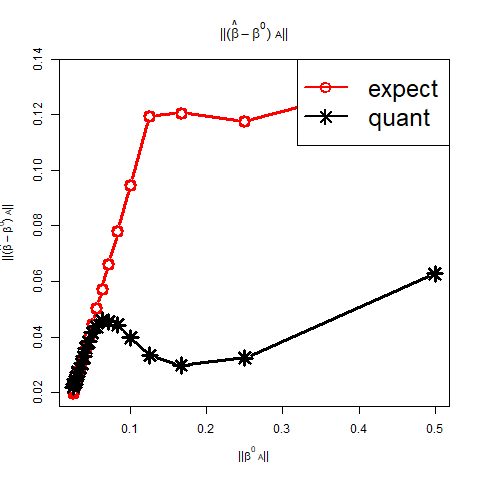} & \\
		{\small (c) Evolution of $L1A\equiv \|\big(\overset{\vee}{\eb}_{n}-\ebo \big)_{\cal A} \|_1$.} & \\
	\end{tabular}
	\caption{\small Study of false zeros and of $\|\big(\overset{\vee}{\eb}_{n}-\ebo \big)_{\cal A} \|_1$ with respect to values of $\|\eb^0_{\cal A}\|_1$ values, when ${\cal A}=\{1\}$ and $\beta^0_1=1/(2j)$, $j \in \{1, \cdots , 20\}$, $K=25$, $w=[\sqrt{K}]$.	}
	\label{fig_evol_beta0}
\end{figure}

\begin{figure}[h!] 
	\begin{tabular}{cc}
		\includegraphics[width=0.45\linewidth,height=4.5cm]{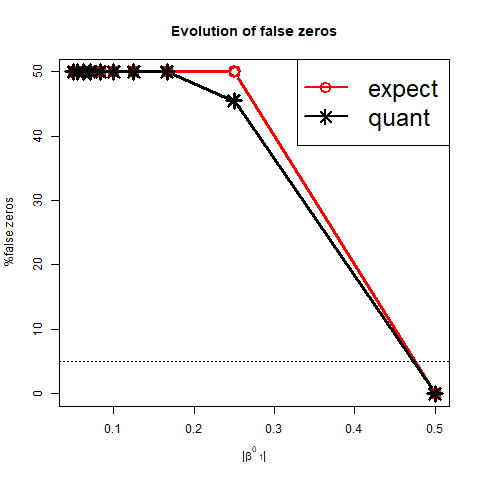} &
		\includegraphics[width=0.45\linewidth,height=4.5cm]{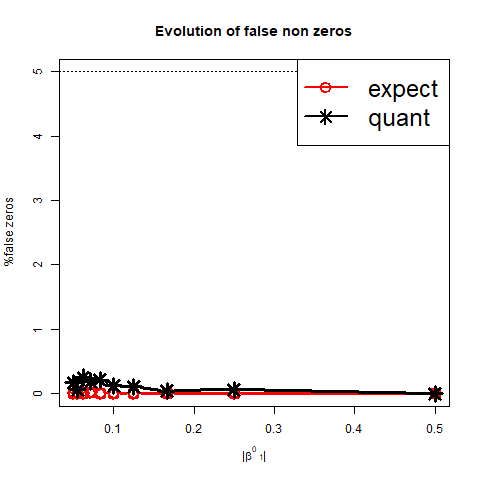} \\
		{\small (a) Evolution of false zeros.} &
		{\small (b) Evolution of false non-zeros.}\\
		\includegraphics[width=0.45\linewidth,height=4.5cm]{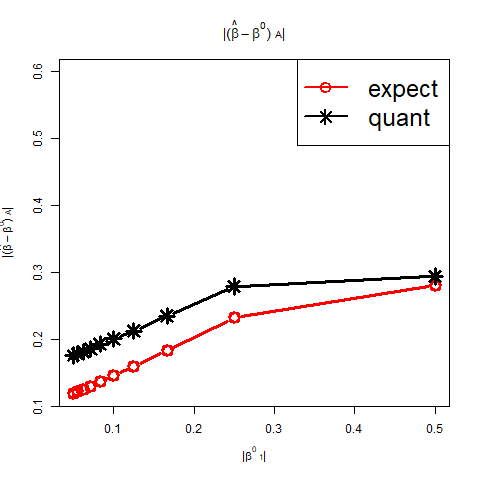} & \\
		{\small (c) Evolution of $L1A=\|\big(\overset{\vee}{\eb}_{n}-\ebo \big)_{\cal A} \|_1$.} & \\
	\end{tabular}
	\caption{\small Study of false zeros, of false non-zeros and of $\|\big(\overset{\vee}{\eb}_{n}-\ebo \big)_{\cal A} \|_1$ with respect to values of  $\|\eb^0_{1}\|_1$, when $p=50$, ${\cal A}=\{1,2\}$ and $\beta^0_1=1/(5j)$, $j \in \{1, \cdots , 10\}$, $\beta^0_2=-2$, $K=25$, $w=[\sqrt{K}]$.	}
	\label{fig_evol_beta2}
\end{figure}

We first study automatic selection of the relevant explanatory variables for aggregated censored adaptive LASSO expectile and quantile estimation methods with respect to the number of non-zero elements of $\ebo$.  For $p=100$, ${\cal A}=\{1, \cdots , 95\}$ with $\beta^0_j=1$ for $j \in \{ 1, 3, 4, \cdots, 95\}$ and $\beta^0_2=-2$, when $n=10^5$,  $K=10$, $w= [\sqrt{K}]$, $\lambda_n=n^{1/2-1/10}$,   $\varepsilon \sim {\cal G}(0,1)$, we obtain that the percentage of false non-zeros is $0$ by the expectile  and  quantile methods. Moreover, the percentage of false zeros is equal to 0 for both methods. In order to compare with the results obtained on all data by the censored adaptive LASSO expectile and quantile estimator $\widehat \eb_n$ we get that: rate of false zeros is equal to 0 by both methods, rate of false non-zeros is  equal to $7.2\%$ and $2.4\%$ by expectile and quantile, respectively. \\
Compare now these results with those obtained for $\beta^0_1=1$, $\beta^0_2=-2$ and $\beta^0_j=0$, for any $  j \in \{3, \cdots , 100\}$.  On all data, the percentage of false zeros is equal to 0, by both methods, and the percentage of false non-zeros is equal to 0.12 and 3.18, by expectile and quantile, respectively. For $\overset{\vee}{\eb}_{n}$ we have $\overset{\vee}{\cal A}_n ={\cal A}=\{1, 2\}$ for all Monte Carlo replications (so the percentage of false zeros is equal to 0). \\
We now study the evolution of automatic variable selection with $\|\ebo\|_1$. First we take ${\cal A}=\{1\}$ with $p=50$, $\beta^0_1=1/(2j)$, $j \in \{1, 2, \cdots, 20\}$, $K=25$, $w=[\sqrt{K}]$ (Figure \ref{fig_evol_beta0}). Note that the false zero and false non-zero rates are 0 or very close to 0. (Figures \ref{fig_evol_beta0}(a) and (b)). The bias of $\overset{\vee}{\eb}_{n, {\cal A}}$  is slightly greater by the quantile method   for values of $\|\eb^0_{\cal A}\|_1$ far from 0, but then it is the same as with the expectile method if  $\|\eb^0_{\cal A}\|_1$ is close to 0 (Figure \ref{fig_evol_beta0}(c)).  In Figure \ref{fig_evol_beta2} we consider ${\cal A}=\{1, 2\}$ with $p=50$ and $\beta^0_1=1/(5j)$, $j \in \{1, \cdots , 10\}$, $\beta^0_2=-2$. We deduce that the evolution of false zero and false non-zero rates  is more or less the same for both estimations   (Figures  \ref{fig_evol_beta2}(a) and \ref{fig_evol_beta2}(b)).  By comparing $L1A \equiv\|\big(\overset{\vee}{\eb}_{n}-\ebo \big)_{\cal A} \|_1$ values when $\beta_2^0=-0$ (Figure \ref{fig_evol_beta0}(c)) and $\beta_2^0=-2$ (Figure \ref{fig_evol_beta2}(c)), we deduce that for $\beta^0_1$ close to 0 when there is another parameter different from 0, then the expectile method gives estimations with a lower bias than quantile, a difference that diminishes as $\beta^0_1$ moves away from 0.

\section{Proofs}
\label{section_proofs}
In this section we present the proofs of the four theorems stated in Section \ref{section_estimators}.\\
\noindent  {\bf Proof of Theorem \ref{th_oracle}}.  \\
\textit{(i)} 
Under assumptions (A1), (A4), (A5)-(A8) together with relation (\ref{eln}), by Theorem 2(i) of \cite{Shows.Lu.Zhang.10} together the fact that  $K=o(n)$, we have for any $k \in \{1, \cdots, K\}$ that  $\PP [\widehat {\cal A}^{(k)}_n=\widehat{\cal A}_n ={\cal A}] {\underset{n \rightarrow \infty}{\longrightarrow}} 1$. Then, taking into account the definition of $\overset{\vee}{\cal A}_n$, we obtain  $\lim_{n \rightarrow \infty} \PP\big[ \overset{\vee}{\cal A}_n= {\cal A} \big]=1$ and $\lim_{n \rightarrow\infty} \PP[ \overset{\vee}{\eb}_{n, {\cal A}^c}=\textbf{0}_{|{\cal A}^c|}]=1$.\\
Corresponding to relation (\ref{eq3}), we have for the $k$-th group of observations:
\begin{align}
	\big(\widehat \eb_n^{(k)} - \ebo\big)_{\cal A}&=-\bigg(\frac{n}{K}\bigg)^{-1} \big(\eS_{\widehat{\cal A}_n^{(k)}}^{(k)}\big)^{-1} \bigg(\sum_{i \in {\cal U}_k} \es_i\bigg) \bigg(1+o_\PP \bigg(\bigg(\frac{n}{K}\bigg)^{-1/2}\bigg)\bigg) \nonumber \\
	& =-\bigg(\frac{n}{K}\bigg)^{-1} \big(\eS_{ {\cal A} }^{(k)}\big)^{-1} \bigg(\sum_{i \in {\cal U}_k} \es_i\bigg) \bigg(1+o_\PP \bigg(\bigg(\frac{n}{K}\bigg)^{-1/2}\bigg)\bigg).
	\label{eq5}
\end{align}
Then, using the definition of $ \overset{\vee}{\eb}_{n}$ and relation (\ref{eq5}), we can write:
\begin{equation*}
	%	\left\{  
	\begin{split} 
		\overset{\vee}{\eb}_{n,{\cal A}}& = \frac{1}{K} \sum^K_{k=1} \widehat \eb_{n,{\cal A}}^{(k)}= \eb^0_{\cal A}+\frac{1}{K} \sum^K_{k=1} \bigg(\frac{n}{K}\bigg)^{-1} \big(\eS_{\cal A}^{(k)}\big)^{-1} \bigg(\sum_{i \in {\cal U}_k} \es_i\bigg)(1+o_\PP((n/K)^{-1/2}))\\
		& =\eb^0_{\cal A}+O_\PP((n/K)^{-1/2}),
	\end{split}
	%	\right. 
\end{equation*}
where we also used the fact that  $\eE[\es_i]=\oo$ and $\Var[\es_i]<\infty$.\\
\textit{(ii)} By Theorem 2 (ii) of  \cite{Shows.Lu.Zhang.10} we have
\[
n^{1/2} \big(\widehat{\eb}_ {n,{\cal A}} -\eb^0_{\cal A}\big) \overset{\cal L} {\underset{n \rightarrow \infty}{\longrightarrow}}  {\cal N} \big(\oo_{|{\cal A}|},\eS^{-1}_{\cal A} \eV_{\cal A} \eS^{-1}_{\cal A} \big).
\]
Then, for the following $|{\cal A}|$-vector, since $K=o(n)$, we have for any $k=1, \cdots, K$, that
\[
\evP_k \equiv (nK^{-1})^{1/2} \big(\eS^{-1}_{\cal A} \eV_{\cal A} \eS^{-1}_{\cal A}\big)^{-1/2} \big( \widehat \eb^{(k)}_n - \ebo\big)_{\cal A} \overset{\cal L} {\underset{n \rightarrow \infty}{\longrightarrow}}  {\cal N} \big(\oo_{|{\cal A}|}, \eI_{|{\cal A}|} \big),
\]
with $\eI_{|{\cal A}|}$ the identity matrix of order $|{\cal A}|$. Taking into account that for $k, k' \in \{1, \cdots, K\}$ the estimators $\widehat \eb^{(k)}_n$ and $\widehat \eb^{(k')}_n$ are independent by the way observation groups ${\cal U}_k$, ${\cal U}_{k'}$ are composed and assumptions (A1), (A2), (A3), we deduce that:
\[
K^{-1/2} \sum^K_{k=1} (\evP_k - \eE[\evP_k])  \overset{\cal L} {\underset{n \rightarrow \infty}{\longrightarrow}}  {\cal N} \big(\oo_{|{\cal A}|}, \eI_{|{\cal A}|} \big).
\]
Then, we obtain then
\[
n^{1/2} \big(\eS^{-1}_{\cal A} \eV_{\cal A} \eS^{-1}_{\cal A}\big)^{-1/2} \big( \overset{\vee}{\eb}_n - \ebo\big)_{\cal A}  \overset{\cal L} {\underset{n \rightarrow \infty}{\longrightarrow}}  {\cal N}_{|{\cal A}|} \big(\oo_{|{\cal A}|}, \eI_{|{\cal A}|} \big)
\]
and the theorem is proved.
\hspace*{\fill}$\blacksquare$ \\

\noindent  {\bf Proof of Theorem \ref{th_oracleq}}.  \\
\textit{(i)} Under assumptions (A1)-(A4), (A5)-(A8) together with relation (\ref{eln}), the proof is identical to that of Theorem \ref{th_oracle}(\textit{i}). \\
\textit{(ii)}
Since $K=o(n)$, then for the estimator $\widehat \eb_{n,{\cal A}}^{(k)}$, taking into account relation (\ref{r123}), we have the relation similar to (\ref{eq5}):
\begin{align}
	\big(\widehat \eb_n^{(k)} - \ebo\big)_{\cal A}&=-\bigg(\frac{n}{K}\bigg)^{-1} \big(\eS_{\widehat{\cal A}_n^{(k)}}^{(k)}\big)^{-1} \bigg(\sum_{i \in {\cal U}_k}  (\es_{i}^{(1)}+\es_{i}^{(2)})\bigg) \bigg(1+o_\PP \bigg(\bigg(\frac{n}{K}\bigg)^{-1/2}\bigg) \bigg)\nonumber \\
	& =-\bigg(\frac{n}{K}\bigg)^{-1} \big(\eS_{{\cal A}}^{(k)}\big)^{-1} \bigg(\sum_{i \in {\cal U}_k} (\es_{i}^{(1)}+\es_{i}^{(2)})\bigg) \bigg(1+o_\PP \bigg(\bigg(\frac{n}{K}\bigg)^{-1/2}\bigg)\bigg).
	\label{eq6}
\end{align}
With respect to the proof of Theorem \ref{th_oracle}(\textit{ii}), in the present case, instead of   $\eV_{\cal A}$, we will have the variance-covariance matrix of $\eW_1+\eW_2$ defined in (\ref{eq_etoile}).   Thus, for the aggregated censored  adaptive LASSO quantile estimator $\overset{\vee}{\eb}_{n}$  we have:
$$n^{1/2} \big(  \overset{\vee}{\eb}_{n,{\cal A}} -\eb^0_{\cal A} \big) \overset{\cal L} {\underset{n \rightarrow \infty}{\longrightarrow}}  {\cal N} \pth{ \textbf{0}_{|{\cal A}|}, \eS_{{\cal A}}^{-1} \Var[\eW_1+\eW_2] \eS_{{\cal A}}^{-1} }.
$$
\hspace*{\fill}$\blacksquare$ \\

\noindent  {\bf Proof of Theorem \ref{th_oracleqJ}}.  \\
\textit{(ii)}
Under assumptions (A1)-(A8) together relation (\ref{eln}), given the proof of Theorem 2.3(ii) of \cite{Tang-Zhou-Wu.12}, under assumptions (A1), (A2), (A4), (A5), (A8),   the estimator  $\widehat{\eb}_n$ of relation (\ref{bn_aquantJ}) is written, with probability converging to 1, as $\ebo+n^{-1/2} \eu$, with $\eu=(\eu_{\cal A}, \oo_{|{\cal A}^c|})$, $\eu_{\cal A} \in \R^{|{\cal A}|}$, $\| \eu_{\cal A}\|_2 \leq c< \infty$.\\
Let also be the vector $\ev=(v_1, \cdots, v_J)$ which comes from the fact that, for all  $j \in \{1, \cdots ,J\}$ we have $n^{1/2}(\widehat b_{n,j} - b^0_{\tau_j})=O_\PP(1)$ (see the proof of Theorem 2.2 of  \cite{Tang-Zhou-Wu.12}), with $b^0_{\tau_j}$ the quantile of order $\tau_j$ of the error $\varepsilon$.\\
Then, $\eu_{\cal A}$ and $\ev$ are minimizers of the following process (see proof of Theorem 2.3 of \cite{Tang-Zhou-Wu.12}), with probability converging to 1:
\[
\frac{1}{2}\eu_{\cal A}^\top \eD_{\cal A} \eu_{\cal A}+\eu_{\cal A}^\top \frac{1}{ \sqrt{n}}\sum^n_{i=1} \es_{i,J}^{(1)}+\eu_{\cal A}^\top \frac{1}{ \sqrt{n}}\sum^n_{i=1} \es_{i,J}^{(2)}+\frac{1}{2}\sum^J_{j=1}\varphi_jv_j^2+\sum^J_{j=1}\frac{v_j}{ \sqrt{n}}\sum^n_{i=1}\frac{\delta_i}{G_0(Y_i)} \big(\e1_{\varepsilon_i < b^0_{\tau_j}-\tau_j}\big),
\] 
from where we get
\begin{equation}
	\label{eq4J}
	\eu_{\cal A}=- \big(\eD_{\cal A}\big)^{-1} \frac{1}{\sqrt{n}}\sum^n_{i=1}\big(\es_{i,J}^{(1)}+\es_{i,J}^{(2)}\big),
\end{equation}
with
\begin{equation}
	\label{eq_etoileJ}
	\left\{  
	\begin{split} 
		n^{-1/2} \sum^n_{i=1} \es_{i,J}^{(1)} \overset{\cal L} {\underset{n \rightarrow \infty}{\longrightarrow}} \eW_{1,J} \sim  {\cal N}_{|{\cal A}|}(\oo_{|{\cal A}|},...)\\
		n^{-1/2}\sum^n_{i=1} \es_{i,J}^{(2)} \overset{\cal L} {\underset{n \rightarrow \infty}{\longrightarrow}} \eW_{2,J} \sim  {\cal N}_{|{\cal A}|}(\oo_{|{\cal A}|},...).
	\end{split}
	\right. 
\end{equation}
The relation similar to (\ref{eq3}) is 
\begin{equation}
	\label{r178}
	\big(\widehat \eb_n - \ebo\big)_{\cal A} =- n^{-1/2}\eD_{\cal A}^{-1} \bigg(n^{-1/2}\sum^n_{i=1} \big(\es_{i,J}^{(1)}+\es_{i,J}^{(2)}\big) \bigg)(1+o_\PP(n^{-1/2})),
\end{equation}
On the observations ${\cal U}_k$ ($k$th group of observations) the censored adaptive LASSO composite quantile estimators corresponding to relations   (\ref{bn_aLAD}) and (\ref{bn_LAD}), are respectively:
\begin{equation}
	\label{bek}
	\left\{  
	\begin{split} 
		\widehat \eb_n^{(k)}& \equiv \argmin_{b_j,\eb \in \R^p} \sum^J_{j=1}  \sum_{i \in {\cal U}_k} \frac{\delta_i}{\widehat G_n(Y_i)} \rho_{\tau_j}\big(\log (Y_i) -b_j- \eX_i^\top \eb\big)+\lambda_{(n_{k-1}+1,n_k)}\sum^p_{j=1} \frac{|\beta_j|}{ |\widetilde \beta^{(k)}_{n,j}|} \\
		\widetilde \eb_n^{(k)}& \equiv \argmin_{b_j, \eb \in \R^p} \sum^J_{j=1}   \sum_{i \in {\cal U}_k} \frac{\delta_i}{\widehat G_n(Y_i)}\rho_{\tau_j}\big(\log (Y_i) -b_j- \eX_i^\top \eb\big).
	\end{split}
	\right. 
\end{equation}
For $\widehat \eb_{n,{\cal A}}^{(k)}$ we have a relation similar to (\ref{eq5}), taking into account relation (\ref{r178}):
\begin{align}
	\big(\widehat \eb_n^{(k)} - \ebo\big)_{\cal A}&=-\bigg(\frac{n}{K}\bigg)^{-1} \big(\eD_{\widehat{\cal A}_n^{(k)}}^{(k)}\big)^{-1} \bigg(\sum_{i \in {\cal U}_k} (\es_{i,J}^{(1)}+\es_{i,J}^{(2)})\bigg) \bigg(1+o_\PP \bigg(\bigg(\frac{n}{K}\bigg)^{-1/2}\bigg) \bigg)\nonumber \\
	& =-\bigg(\frac{n}{K}\bigg)^{-1} \big(\eD_{ {\cal A} }^{(k)}\big)^{-1} \bigg(\sum_{i \in {\cal U}_k}(\es_{i,J}^{(1)}+\es_{i,J}^{(2)})\bigg) \bigg(1+o_\PP \bigg(\bigg(\frac{n}{K}\bigg)^{-1/2}\bigg)\bigg).
	\label{eq6J}
\end{align}
The rest of the proof is similar to that of  Theorem \ref{th_oracle}. In this case, instead of $\eV_{\cal A}$ we will have the variance-covariance matrix of $\eW_{1,J}+\eW_{2,J}$ defined in (\ref{eq_etoile}) and $\eD_{\cal A}$ will replace $\eS_{\cal A}$. So, for the  aggregated censored adaptive LASSO composite quantile estimator $\overset{\vee}{\eb}_{n}$ we have:
$$n^{1/2} \big(  \overset{\vee}{\eb}_{n,{\cal A}} -\eb^0_{\cal A} \big) \overset{\cal L} {\underset{n \rightarrow \infty}{\longrightarrow}}  {\cal N} \pth{ \textbf
	{0}_{|{\cal A}|}, \eD^{-1}_{\cal A} \Var[\eW_{1,J}+\eW_{2,J}] \eD^{-1}_{\cal A}}.
$$
\hspace*{\fill}$\blacksquare$ \\

\noindent  {\bf Proof of Theorem \ref{th_oraclee}}.  \\
\textit{(ii)}
Under assumptions (A1)-(A7), (A9), given Theorem 3(\textit{i}) of  \cite{Ciuperca.24},  consider $\eb=\ebo+n^{-1/2}\eu$, with $\eu=(\eu_{\cal A}, \oo_{|{\cal A}^c|})$, $\| \eu_{\cal A}\| < \infty$. The vector $\eu_{\cal A}$ is the minimizer of (with probability converging to 1, see the proof of Theorem 3 of  \cite{Ciuperca.24}):
\begin{equation*}
	\begin{split} 
		\frac{1}{2} \eE_\varepsilon[h_\tau(\varepsilon)] \eu_{\cal A}^\top \eE_\eX[\eX_{\cal A} \eX_{{\cal A}}^\top] \eu_{\cal A}+\eu_{\cal A}^\top \frac{1}{ \sqrt{n}} \sum^n_{i=1}\bigg(\frac{\delta_i}{G_0(Y_i)}g_\tau(\varepsilon_i) \eX_{{\cal A},i}+\int^B_0 \frac{\ek_{\cal A}(z) }{y(z)}dM^{\cal C}_i(z)\bigg)\\
		+\ell_0 \widehat \eo_{n,{\cal A}} sign(\eb^0_{\cal A}) \eu_{\cal A},
		\qquad \qquad
	\end{split}
\end{equation*}
with $\ell_0 \equiv  \lim_{n \rightarrow \infty} n^{-1/2} \lambda_n$. 
This means
\[
\eu_{\cal A}=-\eE_\varepsilon^{-1}[h_\tau(\varepsilon)]  \eE_\eX[\eX_{\cal A} \eX_{{\cal A}}^\top] ^{-1} \bigg(  n^{-1/2} \sum^n_{i=1} \widetilde \es_i+\ell_0 \widehat \eo_{n,{\cal A}} sign(\eb^0_{\cal A})\bigg).
\]
By assumption (\ref{eln}) the sequence $\lambda_n$ is such that $\ell_0=0$. Then we have a relation similar to (\ref{eq3}):
\begin{equation}
	\label{beA}
\big(\widehat \eb_n - \ebo\big)_{\cal A}=-\eE_\varepsilon^{-1}[h_\tau(\varepsilon)]  \eE_\eX[\eX_{\cal A} \eX_{{\cal A}}^\top] ^{-1}  \frac{1}{ n} \sum^n_{i=1} \widetilde \es_i \big(1+ o_\PP(n^{-1/2})\big).
\end{equation}
On the other hand, taking into account assumptions (A1)-(A4),  by CLT we have:
\[
n^{-1/2} \sum^n_{i=1} \widetilde \es_i  \overset{\cal L} {\underset{n \rightarrow \infty}{\longrightarrow}}  {\cal N}_{|{\cal A}|} (\oo_{|{\cal A}|},\textbf{S}_{1,{\cal A}}+\textbf{S}_{2,{\cal A}}),
\]
i.e. a result similar to (\ref{nb}).\\
For the $k$-th group ${\cal U}_k$ of observations, we have a relation similar to (\ref{eq5}):
\[
\big(\widehat \eb_n^{(k)} - \ebo\big)_{\cal A}=-\bigg(\frac{n}{K}\bigg)^{-1}\eE_\varepsilon^{-1}[h_\tau(\varepsilon)]  \eE_\eX[\eX_{\cal A} \eX_{{\cal A}}^\top] ^{-1}   \sum_{i \in {\cal U}_k} \widetilde \es_i \big(1+ o_\PP((n/K)^{-1/2})\big).
\]
With all these results we can apply a similar approach to the proof of Theorem \ref{th_oracle}, with:
\[
n^{1/2} \big(  \overset{\vee}{\eb}_{n,{\cal A}} -\eb^0_{\cal A} \big) \overset{\cal L} {\underset{n \rightarrow \infty}{\longrightarrow}}  {\cal N}_{|{\cal A}|} \pth{ \textbf{0}_{|{\cal A}|}, \eSS}.
\]
\hspace*{\fill}$\blacksquare$

%\textbf{ORCID}\\
%\textit{Gabriela Ciuperca} https://orcid.org/0000-0001-5467-9774

\end{document}